\input amstex
\documentstyle{amsppt}
\magnification 1200
\vcorrection{-1cm}
\NoBlackBoxes
\input epsf

\rightheadtext{Algorithmic recognition of quasipositive braids}
\topmatter
\title     Algorithmic recognition of quasipositive braids of algebraic length two
\endtitle
\author   S.Yu.~Orevkov
\endauthor
\address  IMT, Universit\'e Toulouse-3, France \endaddress
\address Steklov Math. Institute, Moscow, Russia \endaddress

\abstract
We give an algorithm to decide if a given braid is a product
of two factors which are conjugates of given powers of standard generators
of the braid group. The same problem is solved in a certain class of Garside groups
including Artin-Tits groups of spherical type. The solution is based on the
Garside theory and, especially, on the theory of cyclic sliding developed by
Gebhardt and Gonz\'alez-Meneses. We show that if a braid is of the required form,
then any cycling orbit in its sliding circuit set in the dual Garside
structure contains an element for which this fact is immediately seen
from the left normal form.
\endabstract
\endtopmatter

\def\refBessis         {1}
\def\refBessisCorran   {2} \def\refBC{\refBessisCorran}
\def\refBGGMii         {3} 
\def\refBKLone         {4}
\def\refBKLtwo         {5}
\def\refBrSa           {6}
\def\refCharney        {7}
\def\refCrisp          {8}
\def\refDehornoy       {9}
\def\refDP             {10}       
\def\refDeligne        {11}
\def\refEM             {12}       
\def\refThurstonAndCo  {13}
\def\refGG             {14}
\def\refGarside        {15}
\def\refGebhardt       {16}
\def\refGebGM          {17}
\def\refMurasugi       {18}
\def\refOrevkovTop     {19}
\def\refOrevkovUR      {20}
\def\refOrevkovGAFA    {21}
\def\refOrevkovQPthree {22}
\def\refOrevkovDN      {23}
\def\refOrevkovAa      {24}
\def\refOrW            {25}
\def\refPrasolov       {26}
\def\refRudolph        {27}

\def\sectGarsideOne  {1}
\def\sectGarsideTwo  {2}
\def\sectProofOne    {3}
\def\sectProofTwo    {4}
\def\sectExample     {5} \def\sectExamples{\sectExample}
\def\sectBthree      {6}
\def\sectBlock       {7}

\def\thQP         {1}
\def\thQPone      {1(a)}

\def\thAQP        {2}

\def\defFT           {1.1}

\def\defHomog        {1.2}
\def\defSym          {1.3}
\def\defSqF          {1.4}
\def\defSFLR         {1.5}
\def\defLW           {1.6}
\def\defLCF          {1.7}  \def\defLNF{\defLCF}
\def\convLNF         {1.8}
\def\defSSS          {1.9}
\def\defCycDec       {1.10}
\def\defUSS          {1.11}
\def\defSC           {1.12}
\def\remOrb          {1.13}

\def\lemSBeqRA    {2.1}
\def\lemTrIneq    {2.2}
\def\lemCharneyA  {2.3}
\def\lemLR        {2.4}
\def\defLS        {2.5}
\def\lemCharney   {2.6}
\def\corCharney   {2.7}
\def\lemCharneyR  {2.8}
\def\corCharneyR  {2.9}
\def\lemCDcommut  {2.10}
\def\corCDcommut  {2.11}
\def\lemSCinv     {2.12}
\def\lemSCgcd     {2.13}
\def\defSCG       {2.14}
\def\lemBlackGrey {2.15}
\def\corBlackGrey {2.16}
\def\defTransp    {2.17}
\def\defCat       {2.18}
\def\propCat      {2.19}
\def\corCat       {2.20}

\def\lemSym       {3.1}
\def\lemQPone     {3.2}
\def\lemQi        {3.3}
\def\lemQii       {3.4}
\def\corQii       {3.5}
\def\lemQiii      {3.6}
\def\lemQiv       {3.7}

\def\lemAlcm      {4.1}
\def\lemBrSaSqF   {4.2}
\def\lemBrSa      {4.3}
\def\lemAQPone    {4.4}
\def\lemAG        {4.5}
\def\remAG        {4.6}
\def\lemChain     {4.7}
\def\lemMaxChainA {4.8}
\def\lemMaxChainB {4.9}
\def\lemAQ        {4.10}
\def\lemAQi       {4.11}
\def\lemAQii      {4.12}

\def\lemT         {6.1}
\def\lemTsign     {6.2}
\def\lemTi        {6.3}
\def\remT         {6.4}
\def\propTi       {6.5}
\def\propTiReduc  {\propTi(e)}
\def\propTiBB     {\propTi(f)}

\def\propBlock    {7.1}
\def\corBlock     {7.2}
\def\lemBi        {7.3}
\def\lemBii       {7.4}
\def\lemBiii      {7.5}

\def\eqQPone        {1}
\def\eqAB           {2}
\def\eqSC           {3}
\def\eqAQPone       {4}
\def\eqArtinAB      {5}
\def\eqArtinAx      {6}
\def\eqASC          {7}
\def\eqArtinBx      {8}
\def\eqArtinABx     {9}
\def\eqLNF          {10}
\def\eqQm           {11}
\def\eqQlen         {12}
\def\eqQmin         {13}
\def\eqabab         {14}
\def\eqChain        {15}
\def\eqExOne        {16}
\def\eqTi           {17}
\def\eqTii          {18}
\def\eqMT           {19}
\def\eqTiii         {20}
\def\eqBi           {21}
\def\eqBii          {22}
\def\eqBiii         {23}

\def\RSSS{\operatorname{RSSS}}
\def\SSS {\operatorname{SSS}}
\def\USS {\operatorname{USS}}
\def\SCG {\operatorname{SCG}}
\def\SC  {\operatorname{SC}}
\def\SS  {\operatorname{SS}}
\def\wedgeRev {\wedge^{\!\Lsh}}
\def\veeRev {\vee^{\Lsh}}

\def\Hom {\operatorname{Hom}}
\def\Br {\operatorname{Br}}
\def\Z{\Bbb Z}
\def\iff{\Leftrightarrow}
\def\len{\operatorname{len}}
\def\Qlen{{\len}_{\Cal Q}}
\def\Qmin{\Cal Q_{\min}}
\def\Sign{\operatorname{Sign}}
\def\Null{\operatorname{Null}}
\def\cd{\cdot}

\document

\head Introduction
\endhead

Let $\Br_n$ be the braid group with $n$ strings.
It is generated by $\sigma_1,\dots,\sigma_{n-1}$ (called {\it standard}
or {\it Artin} generators) subject to the relations
$$
   \text{$\sigma_i\sigma_j=\sigma_j\sigma_i$ for $|i-j|>1$;}\qquad
   \text{$\sigma_i\sigma_j\sigma_i=\sigma_j\sigma_i\sigma_j$ for
             $|i-j|=1$}
$$
In this paper we give an algorithm (rather efficient in practice) to decide if a given braid
is the product of two factors which are conjugates of given powers of standard
generators. Since our solution is based on Garside theory, as a by-product we obtain
a solution to a similar problem for a certain class of Garside groups
which includes Artin-Tits groups of spherical type
(we call them in this paper just Artin groups;
note that $\Br_n$ is the Artin group of type $A_{n-1}$).
The main ingredient of our solution is the theory of cyclic sliding
developed by Gebhardt and Gonz\'alez-Meneses in [\refGebGM]. In fact,
we show that if an element $X$ is a product of two conjugates of atom powers,
then its set of sliding circuits $\SC(X)$ contains an element for which this property
is immediately seen from the left normal form. If the Garside structure is symmetric
(which is the case for the dual structures on Artin groups), then
any cycling orbit in $\SC(X)$ contains such an element.

When speaking of Garside groups, we use mostly
the terminology and notation from [\refGebGM].
All  necessary definitions and facts from
the Garside theory are given in \S\sectGarsideOne\ below.
For readers familiar with the Garside
theory, we just say here that by a Garside structure on
a group $G$ we mean a triple $(G,\Cal P,\Delta)$ were $\Cal P$ is the submonoid of
positive elements and $\Delta$ is the Garside element (see details in \S\sectGarsideOne).
The letter length function on $\Cal P$
is denoted by $\|\;\|$ and the set of atoms is denoted by $\Cal A$.

It is convenient also to give the following new definitions.
We say that a Garside structure is {\it symmetric} if, for any two simple elements
$u,v$, one has $(u\prec v)\Leftrightarrow(v\succ u)$. The main example is
the dual Garside structure
on Artin groups introduced by Bessis [\refBessis], see [\refBessis; \S1.2]. In particular, the
Birman-Ko-Lee Garside structure on the braid groups [\refBKLone] is symmetric.
Another example is the braid extension of the complex reflection group $G(e,e,r)$
with the Garside structure introduced in [\refBC].

Following [\refBrSa], we say that $X\in\Cal P$ is {\it square free}
if there do not exist $U,V\in\Cal P$ and $x\in\Cal A$ such that $X=Ux^2V$.
A Garside structure is called {\it square free} if all simple elements are square free.
We say that a Garside structure is
{\it homogeneous} if $\|XY\|=\|X\|+\|Y\|$ for any $X,Y\in\Cal P$, thus, $\|\;\|$
extends up to a unique homomorphism $e:G\to\Z$ such that $e|_{\Cal A}=1$.
Both the standard and the dual Garside structure on Artin groups are square free and
homogeneous.

The conjugacy class of an element $X$ of a group $G$ is denoted by $X^G$.
We use Convention \convLNF\ (see \S\sectGarsideOne\ below)
for the presentation of left (right) normal forms.
Let us give the statements of the main results (the proofs are in \S\sectProofOne\
and in \S\sectProofTwo).

%
%

\proclaim{ Theorem 1 } Let $(G,\Cal P,\delta)$ be a symmetric homogeneous Garside
structure of finite type with set of atoms $\Cal A$. Let $k,l$ be positive integers.
When $k\ge 2$ in Part {\rm(a)} or when $\max(k,l)\ge2$ in Part {\rm(b)}, we suppose in
addition that the Garside structure is square free. Let $X\in G$ and $x,y\in\Cal A$. Then:

\smallskip
{\rm(a)}. $X\in (x^k)^G$
if and only if the left normal form of $X$ is
$$
    \delta^{-n}\cdot A_n\cdot\dots\cdot A_2\cdot A_1\cdot x_1^{k}
    \cdot B_1\cdot B_2\cdot\dots\cdot B_n                     \eqno(\eqQPone)
$$
where $n\ge0$, $x_1\in x^G\cap\Cal A$ and $A_i,B_i$ are simple elements such that
$$
     A_i\delta^{i-1}B_i=\delta^i,  \qquad i=1,\dots,n.         \eqno(\eqAB)
$$

\smallskip
{\rm(b).} $X\in(x^k)^G(y^l)^G$ if and only if either $X\in(x_1^k y_1^l)^G$ or
{\bf any cycling orbit} and  {\bf any decycling orbit}
in the set of sliding circuits $\SC(X)$ {\rm(}see Remark \remOrb\/{\rm)}
contains an element whose left normal form is
$$
    \delta^{-n}\cdot A_n\cdot\dots\cdot A_2\cdot A_1\cdot x_1^{{k}}
    \cdot B_1\cdot B_2\cdot\dots\cdot B_n\cdot y_1^{l}           \eqno(\eqSC)
$$
where $n\ge1$, $x_1\in x^G\cap\Cal A$, $y_1\in y^G\cap\Cal A$, and $A_i,B_i$ are as in Part {\rm(a)}.
\endproclaim

%
%

Thus, under the hypothesis of Theorem 1, we obtain the following
algorithm to decide if a given $X\in G$ belongs to $(x^k)^G(y^l)^G$.
\roster
\item "Step 1."
       Compute $\frak s^i(X)$, $i=1,2,\dots$
       (see Definition \defSC) until $\frak s^i(X)=\frak s^j(X)$
       for some $j<i$. Set $\tilde X=\frak s^i(X)$.
       We have $\tilde X\in\SC(X)$.
\item "Step 2."
       If $\tilde X\in\Cal P$, then check if $\tilde X\in(x_1^k y_1^l)^G$
       for all pairs of atoms $(x_1,y_1)$ in $(x^G)\times(y^G)$ and finish
       the computation.
\item "Step 3."
       Compute $\bold c^i(\tilde X)$, $i=1,2,\dots$ (see Definition \defCycDec)
       until $\bold c^i(\tilde X)=\tilde X$.
       If some of $\bold c^i(\tilde X)$ is of the
       form (\eqSC), then return YES. Otherwise return NO.
\endroster

%
%

\proclaim{ Theorem 2 } Let $(G,\Cal P,\Delta)$ be the standard Garside structure
on an Artin-Tits group of spherical type. Let $k,l$ be positive integers,
$X\in G$ and $x,y\in\Cal A$. Then:

\smallskip
{\rm(a)}. $X\in (x^k)^G$
if and only if the left normal form of $X$ is either $x_1^k$ or
$$
    \Delta^{-n}\cdot A_n\cdot\dots\cdot A_2\cdot A_1\cdot x_1^{k-1}
    \cdot x_1B_1\cdot B_2\cdot\dots\cdot B_n                     \eqno(\eqAQPone)
$$
where $n\ge 1$, $x_1\in x^G\cap\Cal A$ and $A_i,B_i$ are simple elements such that
$$
     A_i\Delta^{i-1}B_i=\Delta^i,  \qquad i=1,\dots,n.         \eqno(\eqArtinAB)
$$
and
$$
     A_1= A'_1 x_1, \qquad A'_1\in\Cal P.                            \eqno(\eqArtinAx)
$$

\if01
\smallskip
{\rm(b).} $X\in(x^k)^G(y^l)^G$ if and only if either $X\in(x_1^k y_1^l)^G$ or
the set of sliding circuits $\SC(X)$
contains an element whose left normal form is
$$
    \Delta^{-n}\cdot y_1A_n\cdot\dots\cdot A_2\cdot A_1\cdot x_1^{k-1}
    \cdot x_1B_1\cdot B_2\cdot\dots\cdot B_n\cdot y_1^{l-1}           \eqno(\eqASC)
$$
where $x_1\in x^G\cap\Cal A$, $y_1\in y^G\cap\Cal A$, and $A_i,B_i$ are simple
elements which satisfy {\rm(\eqArtinAB)}, {\rm(\eqArtinAx)}, and
$$
     B_n= B' x_1, \qquad B'\in\Cal P                            \eqno(\eqArtinBx)
$$
{\rm(}when $n=1$, the expression $x_1B_1\cdot B_2\cdot\dots\cdot B_ny_1$
should be understood as $x_1B_1y_1${\rm)}.
\fi

\smallskip
{\rm(b).} $X\in(x^k)^G(y^l)^G$ if and only if either $X\in(x_1^k y_1^l)^G$ or
the set of sliding circuits $\SC(X)$
contains an element whose left normal form is
$$
    \Delta^{-n}\cdot A_n\cdot\dots\cdot A_1\cdot x_1^{k-1}\cdot 
    x_1B_1\cdot B_2\cdot\dots\cdot B_{n-1}\cdot B_n y_1\cdot y_1^{l-1}    \eqno(\eqASC)
$$
where $n\ge 1$, $x_1\in x^G\cap\Cal A$, $y_1\in y^G\cap\Cal A$, and $A_i,B_i$ are simple
elements which satisfy {\rm(\eqArtinAB)}, {\rm(\eqArtinAx)}, and
$$
     A_n = \tilde y_1 A'_n, \qquad
     \tilde y_1\Delta^n=\Delta^n y_1, \quad
     A'_n\in\Cal P.
                                                                      \eqno(\eqArtinBx)
$$
When $n=1$, the expression $x_1B_1\cdot B_2\cdot\dots\cdot B_ny_1$ in {\rm(\eqASC)}
is understood as $x_1B_1y_1$ and conditions {\rm(\eqArtinAx)} and {\rm(\eqArtinBx)}
should be replaced by
$$
     A_1= \tilde y_1 A''_1 x_1, \qquad
     \tilde y_1\Delta=\Delta y_1, \quad
     A''_1\in\Cal P.                                                 \eqno(\eqArtinABx)
$$
\endproclaim

%
%

\proclaim{ Corollary 3 } Under the hypothesis of Theorem 1 {\rm(}resp. of Theorem 2\/{\rm)},
if $X\in(x^k)^G(y^l)^G$ and $\inf_s X<0$, then
$\ell_s(X) = -2\inf_s X + k+l$ {\rm(}resp. $\ell_s(X) = -2\inf_s X + k+l-2\,${\rm)},
see Definition \defSSS.
\endproclaim


{\bf Remarks.}
(1). In Theorem 1(b) we typed the words ``any cycling/decycling orbit'' in boldface
because this is a very important difference between Theorems 1 and 2.
A computation of a single cycling or decycling orbit is much easier than a
computation of the whole set of sliding circuits.
Moreover, though $\SC(X)$ for a random $X$
is usually not very big, there are examples of reducible (see [\refGebGM; Prop.~9])
and even rigid pseudo-Anosov (see [\refPrasolov])
braids $X\in \Br_n$ of letter length $l=O(n)$ such that $|\SC(X)|$
is exponentially large.
In contrary, the size of a single cycling orbit of a rigid braid is, of course,
bounded by $l$. It seems plausible that the size of any cycling
orbit of any pseudo-Anosov braid is bounded by a polynomial in $n$, $l$.

\smallskip
(2). In applications for real algebraic curves, the standard Garside structure
is more natural than the dual one. So, it would be interesting to prove the analog
of Theorem 2(b) with any cyclic orbit instead of the whole $\SC(X)$.

\smallskip
(3). Theorem 2 extends to any square free homogeneous Garside structures for which
Lemma \lemBrSa, and Lemma \lemAG\ hold (the latter is not needed for Theorem 2(a)).

 \smallskip
(4). It seems plausible that Theorem 1(b) (at least for the braid group)
 remains true if one replaces the words ``any cycling orbit in $\SC(X)$''
 by ``any cycling orbit in $\USS(X)$''.
\smallskip

(5). We say that a braid in $\Br_n$ is {\it quasipositive} if it is a product
of conjugates of standard generators. The {\it quasipositivity problem} (QPP) 
in $\Br_n$ is the algorithmic problem to decide if a given braid is quasipositive or not.
This problem appears very naturally in the study of plane real or complex
algebraic curves (see, e.~g., [\refRudolph], [\refOrevkovTop -- \refOrW]).
It is solved for $n=3$ in [\refOrevkovQPthree] (see \S\sectBthree).

\smallskip
(6). Let $e:\Br_n\to\Z$ be as above,
i.~e., $e\big(\prod_j \sigma_{i_j}^{k_j}\big) = \sum k_j$. Usually, $e(X)$ is called
the {\it algebraic length} of $X$ or the {\it exponent sum} of $X$.
If a braid $X$ is quasipositive, i.~e., if
$X=\prod_{j=1}^k a_j^{-1}\sigma_{i_j}a_j$, then evidently $k=e(X)$.
So, in the case $e(X)<0$ the braid $X$ is never quasipositive;
in the case $e(X)=0$ it is quasipositive if and only if it is trivial (thus
QPP is just the word problem), and if
$e(X)=1$, then QPP is a particular case of the conjugacy problem in $\Br_n$
which is solved by Garside [\refGarside] but in this case the solution is
particularly fast. Indeed, by [\refBKLtwo], ElRifai-Morton's algorithm [\refEM]
gives the result after $\le\|\delta\|\ell(X)$ cyclings where $\ell(X)$ is
the canonical length of $X$ (see Definition \defLCF) and
Theorem 1(a) shows that $\ell(X)/2$ cyclings is enough.
The next case $e(X)=2$ is covered
by Theorem 1(b) or 2(b).

\smallskip
(7). QPP is a particular case of the {\it class product problem} (CPP) -- the algorithmic
problem to decide if a given element of a group belongs to the product of a
given collection of conjugacy classes. CPP in $\Br_n$ for conjugacy classes
of the braids of algebraic singularities also naturally arises in the study
of plane algebraic curves. So, our result is a solution of CPP in $\Br_n$ for
the product of two braids of singularities of type $A_n$.
Since the Artin group
of type $B_n$ is isomorphic to the group of braids with a distinguished string (see
[\refCrisp; Prop.~5.1]), this case is also important for applications to
plane real algebraic curves, especially, when using the method of cubic resolvents (see
[\refOrevkovAa; \S4 and Apdx. A, C]).

\medskip

\definition{ Example } It is shown in [\refOrevkovAa; \S4.4] that the arrangement
of a real pseudoholomorphic quintic curve in $\Bbb{RP}^2$ with respect two lines
shown in [\refOrevkovAa; Fig.~16.12 or Fig.~25.1] is algebraically unrealizable.
The proof is based on the fact that $X\not\in\sigma_1^G(\sigma_1^4)^G$
where $G=\Br_4$ and
$
   X=\Delta^4\big(\sigma_3^2 \sigma_1^{-1}\sigma_2 \sigma_1\sigma_3^2  \sigma_2^{-1}\sigma_1
     \sigma_2^3\sigma_3^2\sigma_2^4\sigma_3^2\sigma_1\sigma_2\sigma_1\big)^{-1}.
$
This fact was proven in [\refOrevkovAa] using a mixture of Burau and Gassner representations.
We have $(\inf_s X,\ell_s(X))=(-6,12)$ for the standard Garside
structure and $(-6,14)$ for the dual one. Thus
the result follows from Corollary 3 in both cases.
\enddefinition

In \S\sectExample\ we give an example which shows the difficulties in
the Garside-theoretical approach to QPP for $e(X)\ge3$. In \S\sectBthree\ we
give an algorithm for QPP in $\Br_3$ and a C program with its implementation.
In \S\sectBlock\ we prove a property of the dual Garside structures which we hope to be
useful for QPP in the general case.

\smallskip

\subhead Acknowledgment \endsubhead
I am grateful to the referee for indicating some mistakes in the first version
of the paper and for many very useful advises.

%
%

\head \S\sectGarsideOne. Elements of Garside theory
     needed for the statement of Theorems 1 and 2
\endhead

Given a group $G$ and $x,y\in G$, we denote $x^y = y^{-1}xy$ and $x^G=\{x^z\mid z\in G\}$.
Garside groups were introduced in [\refDehornoy, \refDP] as a class of groups to which the
technique initiated by Garside [\refGarside] and further developed
in [\refGarside, \refBrSa, \refDeligne, \refCharney, 
\refThurstonAndCo, \refEM, \refBKLone, \refBKLtwo] can be extended.
When speaking of Garside groups, we 
use mostly definitions and notation from [\refGebGM].
For the reader's convenience we give a summary in this section.
A group $G$ is said to be a {\it Garside group} with 
{\it Garside structure} $(G,\Cal P,\Delta)$ if it admits a submonoid
$\Cal P$ satisfying $\Cal P\cap\Cal P^{-1}=\{1\}$, called the
monoid of {\it positive elements}, and a special element $\Delta\in P$
called the {\it Garside element}, such that the following properties hold:

\roster
\item"(G1)"
  The partial order $\preccurlyeq$ defined on $G$ by
  $a\preccurlyeq b \iff a^{-1}b\in\Cal P$ (which is invariant under left
  multiplication by definition) is a lattice order. That is, for every
  $a,b\in G$ there exist a unique least common multiple $a\vee b$
  and a unique greatest common divisor $a\wedge b$ with respect
  to $\preccurlyeq$.
\item"(G2)"
  The set $[1,\Delta]=\{a\in G\,|\,1\preccurlyeq a\preccurlyeq\Delta\}$,
  called the set of {\it simple elements}, generates $G$.
\item"(G3)"
  Conjugation by $\Delta$ preserves $\Cal P$. That is,
  $(X\in\Cal P)\Longrightarrow (X^\Delta\in\Cal P)$.
\item"(G4)"
  For all $X\in\Cal P\setminus\{1\}$, one has:
$$
   \|X\|=\sup\{k\,|\,\exists\,a_1,\dots,a_k\in\Cal P\setminus\{1\}
   \;\text{ such that }\; X=a_1\dots a_k\}<\infty.
$$
\endroster
If $1\preccurlyeq a\preccurlyeq b$, then we say that $a$ is a {\it prefix} of $b$.
We write $a\prec b$ if $a\preccurlyeq b$ and $a\ne b$.
Similarly to $[1,\Delta]$, we denote:
$]1,\Delta]=[1,\Delta]\setminus\{1\}$,
${[1,\Delta[}=[1,\Delta]\setminus\{\Delta\}$,
${]1,\Delta[}={[1,\Delta[}\setminus\{1\}$.
We define the mappings
$$
     \tau:G\to G,\;\tau(X) = X^\Delta, \qquad \text{and}\qquad
     \partial:[1,\Delta]\to[1,\Delta],\;\; \partial A=A^{-1}\Delta.
$$
We call $\partial A$ and $\partial^{-1}A$ the {\it right} and the {\it left
complement of} $A$ respectively.
It is clear that $\partial^2=\tau|_{[1,\Delta]}$ and thus
$\tau([1,\Delta])=[1,\Delta]=\{a\in G\mid \Delta\succcurlyeq a\succcurlyeq 1\}$.

\if01
A monoid $\Cal P$ is called a {\it Garside monoid} if there
exists a Garside structure $(G,\Cal P,\Delta)$.
In this case $\Delta$ is called a Garside element of $\Cal P$.
\fi

\definition{ Definition \defFT }
A Garside structure $(G,\Cal P,\Delta)$ is said to be
{\it of finite type} if the set of simple elements $[1,\Delta]$
is finite. A group $G$ is called a {\it Garside group of
finite type} if it admits a Garside structure of finite type.
\enddefinition

All Garside structures considered in this paper are of finite type.

An element $a\in\Cal P\setminus\{1\}$ is called an
{\it atom} if $a=bc$ with $b,c\in\Cal P$ implies either $a=1$ or $b=1$.
We denote the set of atoms by $\Cal A$. It is clear that if
$X=a_1\dots a_k$, $a_i\in\Cal P$, $k=\|X\|$, then all $a_i$ are atoms.
So, $\Cal A$ generates $\Cal P$ and $\Cal A\subset[1,\Delta]$.

\definition{ Definition \defHomog }
A Garside structure $(G,\Cal P,\Delta)$ is called
{\it homogeneous} if for any $X,Y\in\Cal P$ one has $\|XY\|=\|X\|+\|Y\|$.
In this case we can define a group homomorphism $e:G\to\Z$
such that $e(\Cal A)=\{1\}$ and $e(X)=\|X\|$ for any $X\in\Cal P$.
\enddefinition

Similarly to $\preccurlyeq$ we define the order $\succcurlyeq$ by
$a\succcurlyeq b\iff ab^{-1}\in\Cal P$. It is obvious that
$a\preccurlyeq b$ is equivalent to $a^{-1}\succcurlyeq b^{-1}$. It follows
that $\succcurlyeq$ is also a lattice order and
$\Cal P=\{X\,|\,1\preccurlyeq X\} = \{X\,|\,X\succcurlyeq 1\}$.
We denote the lcm and gcd of $a$ and $b$
with respect to the lattice order $\succcurlyeq$ by $a\veeRev b$
and $a\wedgeRev b$ respectively.

\definition{ Definition \defSym }
A Garside structure is called {\it symmetric}
if for any simple elements $u,v$ one has
$u\preccurlyeq v\iff v\succcurlyeq u$.
\enddefinition

\definition{ Definition \defSqF }  $X\in\Cal P$ is called {\it square free}
if there do not exist $U,V\in\Cal P$ and $x\in\Cal A$ such that $X=Ux^2V$.
A Garside structure is called {\it square free} if all simple elements are square free.
\enddefinition

Till the end of this section we suppose that $(G,\Cal P,\Delta)$ is a Garside structure
with set of atoms $\Cal A$.

\definition{ Definition \defSFLR } Let $A\in\Cal P$.
As in [\refBKLone, \refEM, \refGarside], we define the {\it starting set} $S(A)$
and the {\it finishing set} $F(A)$:
$$
  S(A) = \{x\in\Cal A\mid x\preccurlyeq A\}, \qquad
  F(A) = \{x\in\Cal A\mid A\succcurlyeq x\}.
$$
If, moreover, $A\in[1,\Delta]$, then, following [\refBKLone],
we define the {\it right complementary set} $R(A)$
and the {\it left complementary set} $L(A)$:
$$
  R(A) = \{x\in\Cal A\mid Ax\preccurlyeq\Delta\}, \qquad
  L(A) = \{x\in\Cal A\mid \Delta\succcurlyeq xA\}.
$$
Or, equivalently, $R(A)=S(\partial A)$ and $L(A)=F(\partial^{-1}A)$.
\enddefinition

\definition{ Definition \defLW }
Given two simple elements $A$, $B$, we say that the decomposition
$AB = A\cdot B$ is {\it left weighted} if $A=AB\wedge\Delta$
which is equivalent to $B\wedge\partial A=1$ or to $S(B)\cap R(A)=\varnothing$.
We say that the decomposition
$AB = A\cdot B$ is {\it right weighted} if $B=AB\wedgeRev\Delta$
which is equivalent to $A\wedgeRev\partial^{-1} B=1$ or to $F(A)\cap L(B)=\varnothing$.
\enddefinition

In particular, for any simple element $A$, the decompositions
$A\cdot 1$ and $\Delta\cdot A$ are left weighted
whereas $A\cdot\Delta$ and $1\cdot A$ are not.

\definition{ Definition \defLNF }
Given $X\in G$, we say that a decomposition
$$
    X=\Delta^p\cdot A_1\cdot A_2\cdot\dots\cdot A_r                \eqno(\eqLNF)
$$
is the {\it left normal form}
of $X$ if $A_i\in{]1,\Delta[}$ for $i=1,\dots,r$ and
$A_i\cdot A_{i+1}$ is left weighted for $i=1,\dots,r-1$.
In this case we define the {\it infimum}, the {\it canonical length} and the
{\it supremum} of $X$ respectively by
$\inf X=p$, $\ell(X) = r$, $\sup X=p+r$. In [\refGarside], $\inf X$ is called
the {\it power} of $X$.
\enddefinition

\definition{ Convention \convLNF }
Given $A,B\in G$, we use both notations $AB$ and $A\cdot B$ for the product in $G$.
However, if a mixed notation is used (e.~g., $X=AB\cdot C\cdot D$)
and we say that this decomposition is left/right weighted or in left/right
normal form, then we mean that the dots separate
simple elements and each consecutive pair of these simple elements is left/right
weighted. If $x$ is an atom and an expression $x^k$ appears in a left/right
normal form (as in Theorems 1 and 2), then it stands for $x\cdot\dots\cdot x$ ($k$ times)
and, of course, $A\cdot x^k\cdot B$ means $A\cdot B$ when $k=0$.
\enddefinition

\definition{ Definition \defSSS } Let $X\in G$. The
{\it summit infimum}, the
{\it summit supremum}, and the
{\it summit length} of $X$ are defined as
$\inf_s X = \max\{\inf Y\mid Y\in X^G\}$,
$\sup_s X = \min\{\sup Y\mid Y\in X^G\}$,
$\ell_s(X) = \min\{\ell(Y)\mid Y\in X^G\}$.
The {\it super summit set} of $X$
is $\SSS(X)=\{Y\in X^G\mid \ell(Y)=\ell_s(X)\}$.
It is shown in [\refEM] that $\ell_s(X) = \sup_s X - \inf_s X$ and thus
$$
   \SSS(X)=\{Y\in X^G\mid \inf Y = {\inf}_s X\text{ and }\sup Y = {\sup}_s X\}.
$$
\enddefinition

\definition{ Definition \defCycDec } Let $X\in G$, $\ell(X)>0$,
and let (\eqLNF) be its
left normal form. We define the {\it initial factor} and the {\it final factor} of
$X$ as $\iota(X) = \tau^{-p}(A_1)$ and $\varphi(X)=A_r$. So, we have
$X = \iota(X)\Delta^p A_2\dots A_{r-1}\varphi(X)$ when $r>1$ and we have
$X = \iota(X)\Delta^p = \Delta^p\varphi(X)$ when $r=1$.
We define the {\it cycling} and the {\it decycling} of $X$ as 
$\bold c(X) = X^{\iota(X)}=\Delta^p A_2\dots A_r\,\iota(X)$ and 
$\bold d(X) = \bold c(X^{-1})^{-1} = X^{\varphi(X)^{-1}}
=A_r\Delta^p A_1\dots A_{r-1}$.
\enddefinition

\definition{ Definition \defUSS } Let $X\in G$. The {\it ultra summit set}
of $X$ is
$$
    \USS(X) = \{ Y\in\SSS(X)\mid \bold c^k(Y)=Y\text{ for some } k>0\}.
$$
The {\it restricted super summit set} of $X$ is
$$
    \RSSS(X) = \{ Y\in\SSS(X)\mid \bold c^k(Y)=\bold d^m(Y)=Y
                  \text{ for some } k,m>0\}.
$$
If $Y\in\USS(X)$, we define the {\it cycling orbit} of $Y$ as
$\{\bold c^k(Y)\mid k\ge 0\}$. Similarly, if $\bold d^k(Y)=Y$ for some
$k\!>\!0$, then we define the {\it decycling orbit} of $Y$ as
$\{\bold d^k(Y)\,|\, k\ge 0\}$.
\enddefinition

\definition{ Definition \defSC }
Let $X\in G$ and let (\eqLNF) be its left normal form.
The {\it preferred prefix} of $X$ is
$\frak p(X)=\iota(X)\wedge\partial(\varphi(X))$. In other words, $\frak p(X)$
is the greatest positive $u$ such that $u\preccurlyeq\iota(X)$ and $\varphi(X)u\preccurlyeq\Delta$.
The {\it cyclic sliding} of $X$ is
$\frak s(X)=X^{\frak p(X)}$. The {\it set of sliding circuits} of $X$ is
$$
    \SC(X) = \{Y\in X^G\mid \frak s^m(Y)=Y\text{ for some } m>0\}.
$$
\enddefinition

\definition{ Remark \remOrb } By [\refGebGM; Prop.~2], we have
$\SC(X)\subset\RSSS(X)$ and if $\ell_s(X)>1$, then $\SC(X)=\RSSS(X)$.
Thus, $\SC(X)$ is a disjoint union of cycling orbits as well as a
disjoint union of decycling orbits.
\enddefinition

%
%

\head \S\sectGarsideTwo. Elements of Garside theory used in
                         the proofs of Theorems 1 and 2
\endhead

Let $(G,\Cal P,\Delta)$ be a Garside structure of finite type with set of atoms $\Cal A$.

\proclaim{ Lemma \lemSBeqRA } Let $A\in[1,\Delta]$ and $B=\partial A$, i.~e.,
$AB=\Delta$. Then $S(B)=R(A)$ and $F(A)=L(B)$.
\endproclaim

\demo{ Proof } $x\in S(B)$ $\Leftrightarrow$
($\exists B'\in\Cal P$, $B=xB'$)  $\Leftrightarrow$
($\exists B'\in\Cal P$, $\Delta=AxB'$)  $\Leftrightarrow$
$x\in R(A)$. Thus $S(B)=F(A)$. Symmetrically, $F(A)=L(B)$. \qed
\enddemo

\proclaim{ Lemma \lemTrIneq } {\rm[\refEM; p.~482]}.
For any $X,Y\in G$ one has $\ell(XY)\le\ell(X)+\ell(Y)$. \qed
\endproclaim

\proclaim{ Lemma \lemCharneyA } {\rm[\refCharney; Lemma 2.4]}. Let $X,Y\in\Cal P$
and let $Y_1=\Delta\wedge Y$. Then $\Delta\wedge(XY) = \Delta\wedge(XY_1)$. \qed
\endproclaim

\proclaim{ Lemma  \lemLR}
Suppose that $X=X_1\cdot\dots\cdot X_n$ is right weighted and $Y=Y_1\cdot\dots\cdot Y_m$
is left weighted. If $\Delta\preccurlyeq XY$, then $\Delta\preccurlyeq X_nY_1$.
\endproclaim

\demo{ Proof }
The condition $\Delta\preccurlyeq XY$ can be rewritten as $\Delta\wedge(XY)=\Delta$.
Hence, by Lemma \lemCharneyA, we have $\Delta=\Delta\wedge(XY)=\Delta\wedge(XY_1)$,
i.~e., $\Delta\preccurlyeq XY_1$ and hence $XY_1\succcurlyeq\Delta$.
Then the analog of Lemma \lemCharneyA\ for $\wedgeRev$
yields $\Delta = \Delta\wedgeRev (XY_1) = \Delta\wedgeRev(X_nY_1)$.
\qed
\enddemo

\definition{ Definition \defLS } The {\it local sliding} is the mapping
$\frak{ls}:[1,\Delta]^2\to[1,\Delta]^2$ defined by
$\frak{ls}(u,v) = (us,s^{-1}v)$ where $s=v\wedge\partial u$. Thus, if $(u',v')=\frak{ls}(u,v)$,
then $u'v'=uv$ and $u'\cdot v'$ is left weighted.
\enddefinition

\proclaim{ Lemma \lemCharney }  {\rm[\refCharney; Prop.~3.1].}
Suppose that $X=A_1\cdot A_2\cdot\dots\cdot A_r$
is in left normal form and let $A_0$ be a simple element.
Then the decomposition $A_0X=A_0'\cdot A_1'\cdot\dots\cdot A_r'$ is left weighted
where the $A'_i$'s are defined recursively together with simple elements $t_0,\dots,t_r$ by
 $t_0=A_0$,  $(A'_{i-1},t_i)=\frak{ls}(t_{i-1},A_i)$, $i=1,\dots,r$, and $A'_r=t_r$.
We have  $A'_i\ne\Delta$ for $i>0$ and $A'_i\ne1$ for $i<r$ {\rm(}but
it is possible that $A'_0=\Delta$ or $A'_r=1${\rm)}.

Thus, if we set $s_i=A_i\wedge\partial t_{i-1}$, then we have
$A_i=s_it_i$ and $A_{i-1}'=t_{i-1}s_i$ for $1\le i\le r$, and
the left normal forms of $X$ and $A_0X$ are:
$$
\split
       X   &= \quad s_1t_1\cdot s_2 t_2\cdot\;\;\dots\;\;\cdot s_rt_r\\
      A_0X &= t_0 s_1\cdot t_1s_2\cdot\dots\cdot t_{r-1}s_r\cdot t_r
\endsplit
$$
where the last factor $t_r$ should be removed if it is equal to $1$.
\endproclaim

\proclaim{ Corollary \corCharney } Let the notation be as in Lemma \lemCharney.
Suppose that $A_j$ is an atom for some $j<r$.
Then $\varphi(t_0X)=A_r$.
\endproclaim

\demo{ Proof } Let $A'_i$, $s_i$ and $t_i$, $1\le i\le r$, be as in Lemma \lemCharney.
Since $A_j=s_j t_j$ is an atom, we have either $s_j=1$ or $t_j=1$.

If $s_j=1$, then $t_j=A_j$. Since $A_j\cdot A_{j+1}=t_j\cdot s_{j+1}t_{j+1}$ is left weighted,
it follows that $s_{j+1}=1$, and we obtain by induction that
$A'_i=A_i$ for $i\ge j$, hence $\varphi(t_0X)=A'_r=A_r$.

If $t_j=1$, then $A'_j=s_{j+1}$, hence $t_{j+1}=1$ because otherwise
$A'_j\cdot A'_{j+1}=s_{j+1}\cdot t_{j+1}s_{j+2}$ would not be left weighted.
Hence $A'_j=A_{j+1}$ and
we obtain by induction $A'_i=A_{i+1}$, $j\le i<r$, and $A'_r=1$.
Thus  $\varphi(t_0X)=A'_{r-1}=A_r$.
\qed
\enddemo

Informally speaking, Lemma \lemCharney\ means that if a product of elements is left weighted
everywhere except the first pair of elements, then it can be put into left normal form in
one passage from the left to the right: first we make left weighted the leftmost pair of elements,
then the next pair, and so on.
Similarly, the next lemma shows that if a product of elements is left weighted
everywhere except the last pair of elements, then it can be put into the left normal form in
one passage from the right to the left.

\proclaim{ Lemma \lemCharneyR } {\rm[\refCharney; Prop.~3.3].}
Suppose that $X=A_1\cdot A_2\cdot\dots\cdot A_r$
is in left normal form and let $A_{r+1}$ be a simple element.
Then the decomposition $XA_{r+1}=A''_1\cdot\dots\cdot A''_{r+1}$is left weighted  where
the $A''_i$'s are defined recursively together with simple elements $A'_1,\dots,A'_r$ by
$A'_{r+1}=A_{r+1}$, $(A'_i,A''_{i+1})=\frak{ls}(A_i,A'_{i+1})$, $i=r,\dots,1$, $A''_1=A'_1$.
We have  $A''_i\ne\Delta$ for $i>1$ and $A''_i\ne1$ for $i\le r$ {\rm(}but
it is possible that $A''_1=\Delta$ or $A''_{r+1}=1${\rm)}.

Thus one has 
$$
\split
  XA_{r+1} &= (A_1\cd A_2\cd A_3\cd\dots\cd A_{r-2}\cd A _{r-1}\cd A  _r) \;A  _{r+1}\\
           &= (A_1\cd A_2\cd A_3\cd\dots\cd A_{r-2}\cd A _{r-1})  (A' _r\cd A''_{r+1})\\
           &= (A_1\cd A_2\cd A_3\cd\dots\cd A_{r-2})  (A'_{r-1}\cd A''_r\cd A''_{r+1})\\
           &\qquad\qquad\qquad\quad\dots\dots\dots\\
           &= \;\;A_1\, (A'_2\cd A''_3\cd\dots\cd A''_{r-2}\cd A''_{r-1}\cd A''_r\cd A''_{r+1})\\
           &= (A''_1\cd A''_2\cd A''_3\cd\dots\cd A''_{r-2}\cd A''_{r-1}\cd A''_r\cd A''_{r+1})
\endsplit
$$
where all the products in the parentheses are left weighted. In particular, the left normal
form of $XA_{r+1}$ is the last line with the factor $A''_{r+1}$ removed if it is equal to $1$.
\qed
\endproclaim

\proclaim{ Corollary \corCharneyR } Let {\rm(\eqLNF)} be the left normal form
of $X$. Let $\tilde A_1=\tau^{-p}(A_1)$.

\smallskip
{\rm(a)}.
Suppose that $2\le j\le r$. Let $Y=A_j\dots A_r\,\tilde A_1$ and let
$A''_j\cdot\dots\cdot A''_r\cdot\tilde A''_1$ be the left weighted
decomposition of $Y$.
Then $\frak s(X) = \Delta^p A''_1 A_2\dots A_{j-1} A''_j\dots A''_r$
where $A''_1=\tau^p(\tilde A''_1)$.

\smallskip
{\rm(b)}.
Suppose that $3\le j\le r$ and $A_{j-1}=BC$ where $B,C\in\Cal P$ and
$C\cdot A_j$ is left weighted. Let $Y=CA_j\dots A_r\,\tilde A_1$ and let
$C''\cdot A''_j\cdot\dots\cdot A''_r\cdot\tilde A''_1$ be the left weighted
decomposition of $Y$.
Then $\frak s(X) = \Delta^p A''_1 A_2\dots A_{j-2} BC''A''_j\dots A''_r$
where $A''_1=\tau^p(\tilde A''_1)$.
\qed\endproclaim

\proclaim{ Lemma \lemCDcommut } \rm{ (See [\refGebGM; Lemma 4]).}
If $\ell(X)\ge 2$ and either $\ell(\bold c(\bold d(X))=\ell(X)$ or
$\ell(\bold d(\bold c(X))=\ell(X)$,
then $\bold c(\bold d(X))=\bold d(\bold c(X))=\frak s(X)$. \qed
\endproclaim

\proclaim{ Corollary \corCDcommut } If $\ell(X)\ge 2$ and
$X\in\SC(X)$, then $\bold c(\bold d(X))=\bold d(\bold c(X))=\frak s(X)$.
\endproclaim

\proclaim{ Lemma \lemSCinv }
$\SC(X)$ is invariant under $\tau$, $\bold c$, and $\bold d$.
\endproclaim

\demo{ Proof }
If $\ell_s(X)=1$, then the statement is evident. If $\ell_s(X)\ge2$, then 
it follows from the fact that $\SC(X)=\RSSS(X)$ (see Remark \remOrb)
combined with Corollary \corCDcommut. \qed\enddemo

\proclaim{ Lemma \lemSCgcd } {\rm[\refGebGM; Prop.~7]}.
Let $X\in G$ and let $s,t$ be elements of $G$ such that
$X^s\in\SC(X)$ and $X^t\in\SC(X)$. Then $X^{s\wedge t}\in\SC(X)$.
\endproclaim

\definition{ Definition \defSCG }
Let $X\in G$ and $s\in\Cal P\setminus\{1\}$. We say that
$s$ is an {\it $\SC$-minimal conjugator} for $X$
if $X^s\in\SC(X)$ and $X^t\not\in\SC(X)$ for any $t$ such that $1\prec t\prec s$.
Since $Y\in\SC(X)$ $\Rightarrow$ $Y^\Delta\in\SC(X)$, it follows from Lemma \lemSCgcd\ that all
$\SC$-minimal conjugators for the elements of $\SC(X)$ are simple elements.
We define the {\it sliding circuits graph} $\SCG(X)$ as the directed graph
whose set of vertices is $\SC(X)$ and whose arrows starting at a vertex $Y$
are the $\SC$-minimal conjugators for $Y$. If $s$ is an $\SC$-minimal conjugator 
for $Y$, then the corresponding arrow connects $Y$ to $Y^s$.
\enddefinition

The following statement is an analog of [\refBGGMii; Th.~2.5] for $\SC(X)$
instead of $\USS(X)$.

\proclaim{ Lemma \lemBlackGrey } Let $X\in\SC(X)$ and let $s$ be an $\SC$-minimal
conjugator for $X$. Then one and only one of the following conditions holds:
\roster
\item $\varphi(X)s$ is a simple element.
\item $\varphi(X)\cdot s$ is left weighted. 
\endroster
\endproclaim

\demo{ Proof } Repeat word-by-word the proof of [\refBGGMii; Th.~2.5]
replacing $\USS$ by $\SC$ and
using Lemma \lemSCinv\ and Lemma \lemSCgcd\ instead of [\refBGGMii; Lemma 2.5]
and [\refBGGMii; Th.~1.13] respectively.
\qed\enddemo

\proclaim{ Corollary \corBlackGrey } Let $X\in\SC(X)$
with $\ell(X)>0$ and let $s$ be an $\SC$-minimal
conjugator for $X$. Then $s$ is a prefix of either $\iota(X)$
or $\partial\varphi(X)$, or both.
\endproclaim

\demo{ Proof } Repeat  word-by-word the proof of [\refBGGMii; Cor.~2.7]. \qed\enddemo

Similarly to [\refBGGMii; \S2], we distinguish two kinds of arrows of the graph $\SCG(X)$.
We say that an arrow $s$ starting at $Y$ is {\it black} if $s$ is a prefix of $\iota(Y)$, and
{\it grey} if it is a prefix of $\partial\varphi(Y)$ or, equivalently, if $\varphi(Y)s$
is a simple element. Note that some arrows may be both black and grey.

\definition{ Definition \defTransp }
Let $X\in G$ and $u\in\Cal P$. We define the {\it $\bold c$-transport} of
$u$ at $X$ as $\bold c_X(u) = \iota(X)^{-1}u\,\iota(X^u)$,
thus $\bold c(X^u)=\bold c(X)^{u'}$ for $u'=\bold c_X(u)$.
Similarly we define the {\it $\frak s$-transport} of
$u$ at $X$ as $\frak s_X(u) = \frak p(X)^{-1}u\,\frak p(X^s)$,
thus $\frak s(X^u)=\frak s(X)^{u'}$ for $u'=\frak s_X(u)$, i.~e.,
the following diagrams commute (arrows are conjugations):
$$
   \matrix X \!\!\! \!\!\!& \overset{\iota(X)}\to{{-}\!{-}\!{\longrightarrow}}&
             \!\!\!\!\!\!\!\!\!\!\!\!\!\bold c(X)\\
   u\downarrow&&  \downarrow\bold c_X(u)\\
   X^u \!\!\!\!\!\! & \underset{\iota(X^u)}\to{{-}\!{-}\!{\longrightarrow}}
            &\!\!\!\!\!\!\!\!\!\!\!\!\bold c(X^u)\endmatrix
\qquad\quad
   \matrix X \!\!\! \!\!\!& \overset{\frak p(X)}\to{{-}\!{-}\!{\longrightarrow}}&
             \!\!\!\!\!\!\!\!\!\!\!\!\!\frak s(X)\\
   u\downarrow&&  \downarrow\frak s_X(u)\\
   X^u \!\!\!\!\!\! & \underset{\frak p(X^u)}\to{{-}\!{-}\!{\longrightarrow}}
            &\!\!\!\!\!\!\!\!\!\!\!\!\frak s(X^u)\endmatrix
$$

\if01
We define the mappings $\bold c,\frak s:G\times\Cal P\to G\times\Cal P$ by setting
$\bold c(X,u) = (\bold c(X),\bold c_X(u))$,
$\frak s(X,u) = (\frak s(X),\frak s_X(u))$. Finally, for $n\ge0$,
we define  $\bold c_X^{(n)}(u)$ and $\frak s_X^{(n)}(u)$ -- 
the {\it iterated $\bold c$-
and $\frak s$-transport of $u$ at $X$} by setting
$\bold c_X^n(X,u)=(\bold c^n(X),\bold c_X^{(n)}(u))$ and
$\frak s_X^n(X,u)=(\frak s^n(X),\frak s_X^{(n)}(u))$.
\fi
\enddefinition

It is pointed out in [\refGebGM; p.~98] that $\SC(X)$ can be viewed as a category and
then $\frak s$ becomes a functor which is a category isomorphism. The same is true
for $\bold c$. Let us give precise definitions and statements.

\definition{ Definition \defCat }
For $X\in G$ we define the {\it sliding circuits category}
$\Cal{SC}(X)$. The set of objects is $\SC(X)$. Given $Y,Z\in\SC(X)$, we define
the set of morphisms from $Y$ to $Z$ as $\Hom(Y,Z)=\{u\in\Cal P\mid Y^u=Z\}$.

\enddefinition

\proclaim{ Proposition \propCat }
{\rm(a)}. The mappings $\bold c,\frak s:\SC(X)\to\SC(X)$ and
$$
    \bold c_Y:\Hom(Y,Z)\to\Hom(\bold c(Y),\bold c(Z)),\qquad
    \frak s_Y:\Hom(Y,Z)\to\Hom(\frak s(Y),\frak s(Z))
$$
define functors of $\Cal{SC}(X)$ to itself.

\smallskip
{\rm(b)}. These functors are automorphisms of the category $\Cal{SC}(X)$.
\endproclaim

\demo{ Proof } (a). Follows from the invariance of $\SC(X)$ under $\bold c$ and $\frak s$
(see Lemma \lemSCinv).
.
(b). Follows from [\refGebhardt; Lemma 2.6] and [\refGebGM; Lemma 8].
\qed\enddemo

Since the functors $\bold c|_{\SC(X)}$ and $\frak s|_{\SC(X)}$ are bijective, we may define
their inverses which we denote by $\bold c^{-1}$ and $\frak s^{-1}$.
If $\ell(X)\ge 2$, then we may define the functor $\bold d:\Cal{SC}(X)\to\Cal{SC}(X)$ by setting
$\bold d = \frak s\circ\bold c^{-1}$. By Corollary \corCDcommut, the restriction of this functor
to the set of object $\SC(X)$ coincides with the decycling operator $\bold d$ defined above.

\medskip\noindent{\bf Remark.}
In fact, we could define the functor $\bold d$ as $\bold c^{-1}\circ\frak s$ as well. We do not know
if these definitions are equivalent or not but any of them is equally good for our purposes
(for the proof of Part (b) of Lemma \lemQiv).
\medskip

Let $M(X)=\{(Y,u)\in\SC(X)\times\Cal P\mid Y^u\in\SC(X)\}$ -- all morphisms of $\Cal SC(X)$.
Let us define $\bold c_*,\frak s_*:M(X)\to M(X)$ by setting
$\bold c_*(X,s)=(\bold c(X),\bold c_X(s))$ and $\frak s_*(X,s)=(\frak s(X),\frak s_X(s))$.
Proposition \propCat\ implies that these mappings are invertible, so we may define $\bold d_*$
as $\frak s_*\circ\bold c_*^{-1}$.

\proclaim{ Corollary \corCat } Let $X\in\SC(X)$ and
let $s$ be an $\SC$-minimal conjugator for $X$.
Let $(X',s')$ be $\bold c_*^m(X,s)$,  $\bold d_*^m(X,s)$, or  $\frak s_*^m(X,s)$, $m\in\Z$.
Then $s'$  is an $\SC$-minimal conjugator for $X'$.

In particular, $\bold c_*^m$, $\bold d_*^m$, and $\frak s_*^m$ define automorphisms of the graph $\SCG(X)$.
\qed\endproclaim

%
%

\head \sectProofOne. Symmetric homogeneous case: proof of Theorem 1
\endhead

Let $(G,\Cal P,\delta)$ be a {\bf symmetric homogeneous} Garside structure
of finite type with set of atoms $\Cal A$.
The following simple observation will be used again and again in this section.

\proclaim{ Lemma \lemSym } Let $x$ be an atom and $A$ a simple element.
If $x\in L(A)$, then there exists $x_1\in x^G\cap\Cal A$ such that
$xA = Ax_1$ and hence $x^k A=Ax_1^k$ for any $k$.

If $x\in R(A)$, then there exists $x_1\in x^G\cap\Cal A$ such that
$A x^k = x_1^k A$ for any $k$.
\endproclaim

\demo{ Proof } Let  $x\in L(A)$. Then $xA$ is a simple element.
Since the Garside structure
is symmetric, we have $A\preccurlyeq xA$, i.~e., $xA = Ax_1$ for some $x_1\in\Cal P$.
Since, moreover, the Garside structure is homogeneous, we have
$\|x_1\| = \|Ax_1\|-\|A\| = \|xA\|-\|A\| = \|x\|=1$, thus $x_1\in\Cal A$.
Since $x_1=x^A$, we have $x_1\in x^G$. The case $x\in R(L)$ is similar. \qed
\enddemo

Note that for any $u\in G$, $k\in\Bbb Z$,  we have $(x^k)^u = (x_1^k)^P$
where $u=\delta^{\inf u}P$
(thus $\inf P=0$) and $x_1 = \tau^{\inf u}(x)\in x^G\cap\Cal A$. Hence,
Part (a) of Theorem \thQP\ is an immediate consequence from the following fact.

\proclaim{ Lemma \lemQPone } Under the hypothesis of Theorem 1, suppose that
$X=(x_1^k)^P$ with $x_1\in x^G\cap\Cal A$ and $\inf P=0$.
Let $P= B_1\cdot\dots\cdot B_n$, $n\ge 1$, be the left normal form of $P$
and let $A_1,\dots,A_n$ be defined by {\rm(\eqAB)}.
Then either {\rm(\eqQPone)} is the left normal form of $X$ or there exist
$x_2\in x^G\cap\Cal A$ and $Q\in\Cal P$ such that $X=(x_2^k)^Q$, $\|Q\|<\|P\|$,
and $\ell(Q)\le\ell(P)$.
\endproclaim

\demo{Proof}
Suppose that such $x_2$ and $Q$ do not exist. Let us show that (\eqQPone) is
left weighted.
We should check that if $C_1$ and $C_2$ are two successive factors in (\eqQPone)
(not including $\delta^{-n}$), then $R(C_1)\cap S(C_2)=\varnothing$.
We consider all possible cases for $(C_1,C_2)$.

\smallskip
Case 1. $(C_1,C_2)=(B_i,B_{i+1})$. Follows from the fact
that $B_1\cdot\dots\cdot B_n$ is the left normal form of $P$.

\smallskip
Case 2. $(C_1,C_2)=(x_1,B_1)$. Suppose that $y\in R(x_1)\cap S(B_1)$.
Since $y\in S(B_1)$,
we have $y\preccurlyeq B_1\preccurlyeq P$. Hence $P=yQ$ with $Q\in\Cal P$,
$\|Q\|<\|P\|$, and $\ell(Q)\le\ell(P)$.
Since $y\in R(x_1)$, we have $x_1\in L(y)$.
By Lemma \lemSym, this implies $x_1 y = y x_2$ for some $x_2\in x^G\cap\Cal A$.
and we obtain $X=P^{-1}x_1^kyQ=P^{-1}yx_2^kQ=Q^{-1}x_2^k Q$. Contradiction.

\smallskip
Case 3. $(C_1,C_2)=(x_1,x_1)$. 
Follows from the condition that the Garside structure is
square free when $k\ge 2$.

\smallskip
Case 4. $(C_1,C_2)=(A_1,x_1)$. Suppose that $R(A_1)\cap S(x_1)\ne\varnothing$.
Since $S(x_1)=\{x_1\}$, this means that $x_1\in R(A_1)$. Hence
$A_1 x_1= x_2 A_1$ for some $x_2\in x^G\cap\Cal A$
by Lemma \lemSym.
Thus, denoting $B_2\dots B_n$ by $Q$, we obtain
$X=Q^{-1}\delta^{-1}A_1x_1^kB_1Q=Q^{-1}\delta^{-1}x_2^kA_1B_1Q=Q^{-1}x_3^k Q$ for
$x_3=\tau(x_2)\in x^G\cap\Cal A$. Evidently, $\|Q\|<\|P\|$, and $\ell(Q)<\ell(P)$.
Contradiction.

\smallskip
Case 5. $(C_1,C_2)=(A_{i+1},A_i)$. Follows from the fact that
$B_i\cdot B_{i+1}$ is left weighted (see, e.~g., [\refBGGMii; Remark 1.8] or
[\refEM; proof of Prop.~4.5]).
\qed\enddemo

The rest of this section is devoted to the proof of Theorem 1(b).
So, let us fix $x,y\in\Cal A$ and $k,l\ge1$.
Let
$$
    \Cal Q_m = \{ P^{-1}x_1^kPy_1^l\mid
            \ell(P)\le m,x_1\in x^G,y_1\in y^G\},                \eqno(\eqQm)
$$
For any $X\in(x^k)^G(y^l)^G$ we set
$$
    \Qlen(X) = \min\{m\mid \Cal Q_m\cap X^G\ne\varnothing\},            \eqno(\eqQlen)
$$
$$
    \Cal Q_{\min}(X)=\Cal Q_n\cap X^G\text{ where }n=\Qlen(X).         \eqno(\eqQmin)
$$
If $\Qlen(X)=0$, then the conclusion of Theorem 2(b) holds by definition of 
$\Qlen(X)$, so we shall consider the case when $\Qlen(X)>0$.

\smallskip
From now on $x_1,x_2,\dots$ and $y_1,y_2,\dots$ will always denote some atoms which
are conjugate to $x$ and $y$ respectively.

\proclaim{ Lemma \lemQi }
If $X\in\Qmin(X)$ and $\Qlen(X)>0$,
then the left normal form of $X$ is as stated in Theorem {\rm1(b)}
with $n=\Qlen(X)$.
\endproclaim

\demo{Proof}
Let $X\in\Qmin(X)$. Then $X=P^{-1}x_1^k Py_1^l$ with
$\ell(P)=n=\Qlen(X)$.
Without loss of generality we may assume that $\inf P=0$
(otherwise we replace $x_1$ by $\tau^{\inf P}(x_1)$)
and $\|P\|$ is the minimal possible among all
presentations of $X$ in this form.
Let $P=B_1\cdot\dots\cdot B_n$ be the left normal form of $P$ and
let $A_1,\dots,A_n$ be defined by (\eqAB).
Then (\eqSC) represents $X$. Let us show that (\eqSC) is left weighted.
By Lemma \lemQPone, the part $\delta^{-n}\cdot A_n\cdot\dots\cdot B_n$ of (\eqSC) is left weighted,
so, it remains to prove that $B_n\cdot y_1$ is left weighted. Suppose that it is not.
Then $y_1\in R(B_n)$ and, by Lemma \lemSym, we obtain
$B_n y_1^l = y_2^l B_n$.
Thus $X$ is conjugate to $B_n\delta^{-n}A_n\dots A_1 x_1^kB_1\dots B_{n-1}y_2^l
=\delta^{-(n-1)}A_{n-1}\dots A_1 x_1^k B_1\dots B_{n-1} y_2^l$
which contradicts the fact that $n=\Qmin(X)$.
\qed\enddemo

\proclaim{ Lemma \lemQii }
If $X\in\Qmin(X)$ and $\Qlen(X)>0$, then $\frak s(X)\in\Qmin(X)$.
\endproclaim

\demo{Proof}
By Lemma \lemQi, we may assume that the left normal form of $X$ is (\eqSC) with
$n=\Qlen(X)$. Let $A=A_{n-1}\dots A_1$ and $B=B_1\dots B_{n-1}$.
Let $u=\frak p(X)$ (see Definition \defSC).
Then we have $A_n=\tau^{-n}(u)A'_n$, 
$A'_n\in\Cal P$, and
$y_1u\preccurlyeq\delta$. In particular, we have $y_1\in L(u)$,
hence Lemma \lemSym\ implies $y_1^l u=uy_2^l$.
By (\eqAB) we have also $\tau^{-n}(u)A'_n\delta^{n-1} B_n=\delta^n$ which is equivalent to
$\tau^{n-1}(A'_n)B_nu=\delta$. Thus $B_nu$ is a simple element and we obtain
$\frak s(X) = \delta^{-n}A'_nA x_1^k B B_ny_1^l u
 = \delta^{-n}A'_nA x_1^k B B_n u y_2^l = P^{-1}x_1^k P y_2^l$ where
$P$ is a product of $n$ simple elements: $P=B_1\cdot\dots\cdot B_{n-1}\cdot B_nu$.
Hence $\ell(P)\le n$ and we obtain $\frak s(X)\in\Cal Q_n=\Qmin(X)$.
\qed\enddemo

\proclaim{ Corollary \corQii } If $X\in(x^k)^G(y^l)^G$, $\Qlen(X)>0$,
then $\SC(X)\cap\Qmin(X)\ne\varnothing$. \qed
\endproclaim

Thus, $\SC(X)$ contains at least one element of the desired form if $\Qlen(X)>0$.

\proclaim{ Lemma \lemQiii}
Let $X\in\SC(X)\cap\Qmin(X)$, $\Qlen(X)>0$,
and let $s$ be an $\SC$-minimal conjugator for $X$. Then:

\smallskip\noindent
{\rm(a)}. If $\varphi(X)s\prec\delta$, i.~e., if the arrow $X\overset s\to\to X^s$ is grey, then
     either $X^s$ or $\bold c(X^s)$ is in $\Qmin(X)$.

\smallskip\noindent
{\rm(b)}. If $\varphi(X)\cdot s$ is left weighted, i.~e., if the arrow $X\overset s\to\to X^s$ is black,
then $\bold d(X^s)=X$.
\endproclaim

\demo{ Proof }  Let $X=P^{-1} x_1^k P y_1^l$ with $P\in\Cal P$,
$\ell(P)=n=\Qlen(X)$. We have $\ell(X)=k+l+2n$ by Lemma \lemQi.

\smallskip
(a). Since $\varphi(X)s=y_1 s\preccurlyeq\delta$, we have $y_1^l s= s y_2^l$ by Lemma \lemSym.
Hence $X^s = X_0 y_2^l$ where $X_0=(Ps)^{-1} x_1^k (Ps)$.
It follows from Lemma \lemQPone\ that $\ell(X_0)=2m+k$ and
$X^s\in\Cal Q_m$ with $m\le\ell(Ps)\le n+1$.
Since $n=\Qlen(X)$, it follows that $m\ge n$ and if $m=n$, then $X^s\in\Qmin(X)$
and we are done. So, we suppose that $m=n+1$. Then $X_0=\delta^{-(n+1)}X_1$
where $X_1\in\Cal P$ and, by the ``right-to-left version'' of Lemma \lemQPone,
the right normal form of
$X_1$ is $A_{n+1}\cdot\dots\cdot A_1\cdot x_2^k\cdot B_1\cdot\dots\cdot B_{n+1}$
with $A_i,B_i$ satisfying (\eqAB) for $i=1,\dots,n+1$.
Since $X^s\in\SC(X)$, we have $\inf X^s=\inf X=n$ which implies that $\delta\prec X_1y_2^l$.
Since $y_2^l=y_2\cdot\dots\cdot y_2$ is the left normal form of $y_2^l$,
it follows from Lemma \lemLR\ that $\delta\preccurlyeq B_{n+1}y_2$.
Since $\|y_2\|=1$ and $\|B_{n+1}\|<\|\delta\|$, this yields $B_{n+1}y_2=\delta$.
This fact combined with $A_{n+1}\delta^n B_{n+1}=\delta^{n+1}$
implies $A_{n+1}=\tau^{-(n+1)}(y_2)$, thus
$$
\split
   X^s &= \delta^{-(n+1)}\cdot \tau^{-(n+1)}(y_2)\cdot A_n\cdot\dots\cdot A_1\cdot
    x_2^k\cdot B_1\cdot\dots\cdot B_n
    \cdot\delta\cdot y_2^{l-1}
\\
     &= \delta^{-n} \cdot\tau^{-n}(y_2)\cdot\tau\big(A_n\cdot\dots\cdot
    A_1\cdot x_2^k\cdot B_1\cdot\dots\cdot B_n\big)
    \cdot y_2^{l-1}
\\
     &= \delta^{-n} \cdot\tau^{-n}(y_2)\cdot A'_n\cdot\dots\cdot
    A'_1\cdot x_3^k\cdot B'_1\cdot\dots\cdot B'_n
    \cdot y_2^{l-1}
\endsplit
$$
where $A'_n\cdot\dots\cdot A'_1\cdot x_3^k\cdot B'_1\cdot\dots\cdot B'_n$ is
the left normal form of $\tau(A_n\dots A_1 x_2^k B_1\dots B_n)$.
The number of simple factors in this decomposition of $X^s$
is equal to $k+l+2n=\ell(X^s)$. Hence, by Lemma \lemCharney,
we have $\iota(X^s)=y_2 t$ with $t\preccurlyeq\tau^n(A'_n)$. 
Then we have $y_2 t=t y_3$ by Lemma \lemSym.
Since $\tau^{-n}(t)\preccurlyeq A'_n$, we have also $A'_n=\tau^{-n}(t)u$
where $u$ is a simple element.
Hence, we obtain
$$
\split
   \bold c(X^s) &= \delta^{-n} u A'_{n-1}\dots A'_1 x_3^k B'_1\dots B'_n y_2^l t
\\
                &= \delta^{-n}\cdot u\cdot A'_{n-1}\cdot\dots\cdot A'_1\cdot x_3^k\cdot
                   B'_1\cdot\dots\cdot B'_{n-1}\cdot
                   B'_n t\cdot y_3^l
\endsplit
$$
Since $u\cdot\delta^{n-1}\cdot B'_n t = \delta^n$, we conclude that $\bold c(X^s)\in\Qmin(X)$.

\smallskip
(b). Let the left normal form of $X$ be as in (\eqSC).
We have $1\prec s\preccurlyeq st =\iota(X)=\tau^n(A_n)$. Hence
$$
  X^s = \delta^{-n}\cdot \tau^{-n}(t)\,\big( A_{n-1}\cdot\dots\cdot A_1\cdot
        x_1^k\cdot B_1\cdot\dots\cdot B_n \cdot y_1^l\cdot s\big).
$$
Since the tail of this decomposition starting with $A_{n-1}$ is left weighted, we have
$\varphi(X^s)=s$ by Corollary \corCharney, hence $\bold d(X^s) = X$.
\qed
\enddemo

\proclaim{ Lemma \lemQiv }
{\rm(a)}. Let $X\in\SC(X)$, $\Qlen(X)>0$,
and let $s$ be an $\SC$-minimal conjugator for $X$.
Suppose that the cycling orbit of $X$ contains an element of $\Qmin(X)$.
Then the cycling orbit of $X^s$ also contains an element of $\Qmin(X)$.

\smallskip\noindent
{\rm(b)}. The same statement for the decycling orbits.
\endproclaim

\demo{ Proof } (a). Let $Y=\bold c^m(X)$ be the element of the $\bold c$-orbit of $X$
which belongs to $\Qmin(X)$. Let $(Y,t)=\bold c^m_*(X,s)$ (see the end of \S\sectGarsideTwo).
By Corollary \corCat, $t$ is an $\SC$-minimal conjugator for $Y$, i.~e.,
$Y\overset t\to\to Y^t$ is an arrow of the graph $\SCG(X)$.

By Corollary \corBlackGrey, any arrow of
$\SCG(X)$ is either grey or black or both grey and black. Hence, by Lemma \lemQiii\ applied
to $Y$ and $t$, one of
$Y^t$, $\bold c(Y^t)$, or $\bold d(Y^t)$ is in $\Qmin(X)$. In the former two cases we are done.
In the latter case it suffices to note that if  $\bold d(Y^t)\in\Qmin(X)$, then
$Z=\frak s^{-1}(\bold d(Y^t))\in\Qmin(X)$ by Lemma \lemQii\ (as in the end of \S\sectGarsideTwo,
here $\frak s^{-1}$ stands for the inverse of $\frak s|_{\SC(X)}$) and 
$Z=\frak s^{-1}(\bold d(\bold c^m(X^s)))=\bold c^{m-1}(X^s)$ by Corollary \corCDcommut, 
thus $Z$ is an element of the cycling orbit of $X^s$ belonging to $\Qmin(X)$.

\smallskip
(b). The same proof but with $\bold c$ and $\bold d$ exchanged.
\qed\enddemo

Theorem 1(b) follows immediately from Lemma \lemQi, Corollary \corQii, and Lemma \lemQiv\
combined with the fact that the graph $\SCG(X)$ is connected (see [\refGebGM; Cor.~10]).

%
%

\head \S\sectProofTwo. Artin groups: proof of Theorem 2
\endhead

Let $(G,\Cal P,\Delta)$ be the standard Garside structure on an Artin-Tits group
of spherical type. This is the case studied in details in [\refBrSa, \refDeligne].
We recall that $G=\langle\Cal A\mid\Cal R\rangle$ where $\Cal A$ can be considered
as the set of vertices of a Coxeter graph
(one of $A_n$, $B_n$, $D_n$, $E_6$, $E_7$, $E_8$, $F_4$, $G_2$, $H_3$, $H_4$, $I_2(p)$)
and $\Cal R=\{R_{ab}\mid a,b\in\Cal A\}$ where
$R_{ab}$ is the relation $\langle a\,b\rangle^{m_{ab}} = \langle b\,a\rangle^{m_{ab}}$.
The notation  $\langle a\,b\rangle^m$ means
$$
     \langle a\,b\rangle^m = \underset\text{$m$ letters}\to{\underbrace{abab\dots}}
      = \cases (ab)^{m/2}, &\text{$m$ is even,}\\
          (ab)^{(m-1)/2}a, &\text{$m$ is odd.}\endcases                  \eqno(\eqabab)
$$
The matrix $(m_{ab})$ is encoded by the Coxeter graph in the usual way.
The set of atoms of the standard Garside structure is $\Cal A$, and $\Cal P$ is the set of
products of atoms.

\proclaim{ Lemma \lemAlcm } {\rm(Follows from [\refBrSa; Lemma 3.3])}.
 $a\vee b = \langle a\,b\rangle^{m_{ab}} = \langle b\,a\rangle^{m_{ab}}$
for  $a,b\in\Cal A$. \qed
\endproclaim

\proclaim{ Lemma \lemBrSaSqF } {\rm[\refBrSa; Lemma 5.4]}. Let $X\in\Cal P$. Then
$X$ is simple if and only if it is square free. \qed
\endproclaim

In our notation, Lemma 3.4 from [\refBrSa] can be reformulated as follows.

\proclaim{ Lemma \lemBrSa } Let $W$ be a simple element of $G$.
Then $S(W)=\Cal A\setminus L(W)$
and $F(W)=\Cal A\setminus R(W)$.
\qed\endproclaim

\noindent
{\bf Remark. } The statement of Lemma \lemBrSa\ is wrong for
the dual Garside structures on the braid groups.
\medskip

The proof of Theorem 2(a) is very similar to that of Theorem 1(a).
It is an immediate consequence of the following fact.

\proclaim{ Lemma \lemAQPone } Under the hypothesis of Theorem \thAQP, suppose that
$X=(x_1^k)^P$ with $x_1\in x^G\cap\Cal A$, $\inf P=0$.
Let $P=B_1\cdot\dots\cdot B_n$, $n\ge 1$, be the left normal form of $P$ and
let $A_1,\dots,A_n$ be defined by {\rm(\eqArtinAB)}.
Then either {\rm(\eqAQPone)} is the left normal form of $X$ and {\rm(\eqArtinAx)} holds,
or there exist
$x_2\in x^G\cap\Cal A$ and $Q\in\Cal P$ such that $X=(x_2^k)^Q$, $\|Q\|<\|P\|$,
and $\ell(Q)\le\ell(P)$.
\endproclaim

\demo{Proof}
Suppose that such $x_2$ and $Q$ do not exist and let us show that $x_1B_1$ is
a simple element, (\eqArtinAx) holds, and (\eqAQPone) is left weighted. Indeed:

\smallskip
Suppose that $x_1B_1$ is not a simple element, i.~e., $x_1\not\in L(B_1)$. By Lemma
\lemSBeqRA\ and Lemma \lemBrSa, this implies $x_1\in S(B_1)=R(A_1)$. Hence $B_1=x_1B'_1$
and we obtain $X=(x_1^k)^Q$ with $Q=B'_1B_2\dots B_n$, $\|Q\|<\|P\|$, and
$\ell(Q)\le\ell(P)$. Contradiction.

\smallskip
Since $x_1\in L(B_1)$, Lemma \lemSBeqRA\  implies that
$x_1\in F(A)$, thus (\eqArtinAx) holds.

\smallskip
Let us show that (\eqAQPone) is left weighted.
We should check that if $C_1$ and $C_2$ are two successive factors in (\eqAQPone)
(not including $\Delta^{-n}$), then $R(C_1)\cap S(C_2)=\varnothing$.
We consider all possible cases for $(C_1,C_2)$.

\smallskip
Case 1. $(C_1,C_2)=(B_i,B_{i+1})$, $i\ge 2$.
Follows from the fact that $B_1\cdot\dots\cdot B_n$ is the left normal form of $P$.

\smallskip
Case 2. $(C_1,C_2)=(x_1B_1,B_2)$. 
Follows from the fact that $B_1\cdot B_2$ is left weighted.

\smallskip
Case 3. $(C_1,C_2)=(\varphi(A_1\cdot x_1^{k-1}),x_1B_1)$.
By (\eqArtinAx) combined with Lemma \lemBrSaSqF, we have $x_1\not\in R(C_1)$.
So, it is enough to show that $S(x_1B_1)=\{x_1\}$.
Suppose that there exists $x_2\in S(x_1B_1)\setminus\{x_1\}$.
Then we have $x_1\preccurlyeq x_1B_1$ and $x_2\preccurlyeq x_1B_1$, 
hence $x_1\vee x_2\preccurlyeq x_1B_1$. Let $x_1B_1 = (x_1\vee x_2)B$.

It follows from Lemma \lemAlcm\ that
$x_1(x_1\vee x_2) = (x_1\vee x_2)x_i$, $i\in\{1,2\}$.
So, by (\eqArtinAx), we have
$A_1 x_1^k B_1 = A'_1x_1^{k+1}B_1 = A'_1x_1^k(x_1 \vee x_2)B=
 A_1'(x_1\vee x_2)x_i^k B=Ax_i^k B$ where $A=A'_1(x_1\vee x_2)$.
Since $AB = A'_1(x_1\vee x_2)B = A'_1x_1B_1 = A_1B_1=\Delta$,
we obtain a contradiction with the minimality of $\|P\|$.

\smallskip
Case 4. $(C_1,C_2)=(x_1,x_1)$. (when $k\ge 3$).
See Lemma \lemBrSaSqF.

\smallskip
Case 5. $(C_1,C_2)=(A_1,x_1)$ (when $k\ge 2$).
Combine (\eqArtinAx) and Lemma \lemBrSaSqF.

\smallskip
{\it Case 6.} $(C_1,C_2)=(A_{i+1},A_i)$. See Case 5 of the proof of Theorem \thQPone.
\qed\enddemo


In our proof of Theorem 2(b) we use one more particular property of Artin groups
which is not a property of any Garside group.

\proclaim{ Lemma \lemAG } Let $a,b\in\Cal A$ and $A\in[1,\Delta]$. If
$a\preccurlyeq Ab$ and
$a\not\preccurlyeq A$, then $Ab=aA$.
\endproclaim

\demo{ Proof } Combine Lemmas \lemChain(b), \lemMaxChainA, and \lemMaxChainB. \qed
\enddemo

\definition{ Remark \remAG } (1).
Let us denote the Artin group corresponding to a Coxeter graph $\Gamma$
by $\Br(\Gamma)$.
In the case when $G$ is the braid group (i.~e.,  $G=\Br(A_n)$),
Lemma \lemAG\ immediately follows from the interpretation of simple elements as
permutation braids given in [\refEM]. Due to the embedding $\Br(B_n)\to\Br(A_{2n})$
(see [\refCrisp; Prop.~5.1]),
the same arguments work also in the case $G=\Br(B_n)$. 

\smallskip
(2). Lemma \lemAG\ can be reformulated as follows: {\it
if $A\in[1,\Delta]$, $y\in\Cal A$, and $\|y\vee A\|\le\|A\|+1$, then
$y\vee A\succcurlyeq A$.}
This statement is no longer true if one omits the
condition $\|y\vee A\|\le\|A\|+1$, Indeed, let $G=\Br_4$, $A=\sigma_2\sigma_1\sigma_3$,
and $y=\sigma_1$. Then we have
$y\vee A=\sigma_2\sigma_1\sigma_3\sigma_2\sigma_3\not\succcurlyeq A$.
\enddefinition

In Lemmas \lemChain\ -- \lemMaxChainB\ below,
we use the divisibility theory for Artin groups developed
by Brieskorn and Saito in [\refBrSa; \S3].
Let us recall some notions and facts from [\refBrSa]. Let $\Cal A^*$ be the free monoid freely
generated by $\Cal A$ (the set of all words in the alphabet $\Cal A$).
Let $a,b\in\Cal A$. We say that a word $C\in\Cal A^*$ is an {\it elementary or
primitive chain from $a$ to $b$} and we write $a\overset C\to\to b$
if there exist $c\in\Cal A\setminus\{a\}$ and $j$, $0<j<m_{ac}$,
such that $C=\langle c\,a\rangle^j$, and $b$ is the last letter of
$\langle c\,a\rangle^{j+1}$, thus $aCb=\langle a\,c\rangle^{j+2}$
(the notation $\langle\dots\rangle^j$ is introduced by (\eqabab)).
The chain $C$ is called {\it primitive} when $m_{ac}=2$ and it is called
{\it elementary} when $m_{ac}>2$. We say that $C$ is {\it saturated} if $j=m_{ac}-1$.

A word $W\in\Cal A^*$ is called a {\it chain from $a$ to $b$} 
if there exists a sequence of elementary or primitive chains
$$
  a=a_0\overset{C_1}\to\longrightarrow a_1
       \overset{C_2}\to\longrightarrow\dots
       \overset{C_n}\to\longrightarrow a_n=b.                           \eqno(\eqChain)
$$
such that $W=C_1\dots C_n$. It is {\it saturated} if each of
$C_1,\dots,C_n$ is a saturated.

\proclaim{ Lemma \lemChain }
{\rm(a)}. Let $a\in\Cal A$ and $X\in\Cal P$. Then one and only one
of the following possibilities holds:
\roster
\item"(i)" $a\preccurlyeq X$.
\item"(ii)" $X$ can be represented by a chain from $a$ to $c$ for some atom $c$.
\endroster

\smallskip\noindent
{\rm(b)}. Let $a,b\in\Cal A$ and $X\in\Cal P$. Suppose that $a\preccurlyeq Xb$
and $a\not\preccurlyeq X$.
Then $X$ can be represented by a chain from $a$ to $b$.
\endproclaim

\demo{ Proof }
(a). Follows from [\refBrSa; Lemma 3.2 and Lemma 3.3].

\smallskip
(b). Since $a\not\preccurlyeq X$, it follows from (a) that
$X$ can be represented by a chain $a\to c$ for some atom $c$.
Suppose that $c\ne b$. Then the chain can be extended up to
$a\to c\overset b\to\to c$ which represents $Xb$.
By (a), this contradicts the condition $a\preccurlyeq Xb$. \qed
\enddemo

\proclaim{ Lemma \lemMaxChainA }
Let $a,b\in\Cal A$ and $X\in\Cal P$. Suppose that $X$ is represented by a saturated chain
from $a$ to $b$. Then $aX=Xb$.
\endproclaim

\demo{ Proof } Suppose that $X$ is represented by an elementary or primitive saturated chain.
Then $X=\langle ca\rangle^{m-1}$, $m=m_{ac}$,
hence $aX=a\langle ca\rangle^{m-1}=\langle ac\rangle^m=\langle ca\rangle^m=\langle ca\rangle^{m-1}b=Xb$.
In the general case, if $X=C_1\dots C_n$ is as in (\eqChain), then
$$
  a_0C_1\dots C_n=C_1a_1C_2\dots C_n=C_1C_2a_2C_3\dots C_n=\dots=C_1\dots C_na_n. \qed
$$
\enddemo

\proclaim{ Lemma \lemMaxChainB }
Let $A$ be a simple element of $G$ represented by a chain $W$ from $a$ to $b$.
If $a\preccurlyeq Ab$, then the chain $W$ is saturated.
\endproclaim

\demo{ Proof } Let $W$ be as in (\eqChain). Let $i$ be the minimal index such that the
chain
$$
    a_i\overset{C_{i+1}}\to\longrightarrow a_{i+1}
       \overset{C_{i+2}}\to\longrightarrow\dots
       \overset{C_n}\to\longrightarrow a_n=b
$$
is saturated. If $i=0$, then we are done. Suppose that $i\ge 1$. Then, for some
$c\in\Cal A\setminus\{a_i\}$ and $j\le m_{a_ic}-2$, 
we have $C_i = {\dots c a_i c}$ ($j$ letters), hence
$C_ia_i = {\dots c a_i c a_i}$ ($j+1$ letters), i.~e., $C_ia_i$ is an elementary chain
from $a_{i-1}$ to $c$.

\smallskip
Case 1. $c\preccurlyeq C_{i+1}\dots C_n$. Since $C_i=(\dots ca_ic)\succcurlyeq c$, it follows that
$A=C_1\dots C_n$ is not square free. Since $A$ is simple, this fact contradicts Lemma \lemBrSaSqF.

\smallskip
Case 2. $c\not\preccurlyeq C_{i+1}\dots C_n$. Then, by Lemma \lemChain(a), we have
$C_{i+1}\dots C_n = C'_1\dots C'_p$ where
$
   c\overset{C'_1}\to\longrightarrow\dots
    \overset{C'_p}\to\longrightarrow d
$
is a chain from $c$ to some atom $d$. By Lemma \lemMaxChainA, we have
$C_{i+1}\dots C_n b = a_i C_{i+1}\dots C_n$, thus
$Ab =  C_1\dots C_i a_i C'_1\dots C'_p$ which means that
$$
   a=a_0\overset{C_1}\to\longrightarrow\dots
        \overset{C_{i-1}}\to\longrightarrow a_{i-1}
        \overset{C_i a_i}\to\longrightarrow c
        \overset{C'_1}\to\longrightarrow\dots
        \overset{C'_p}\to\longrightarrow d
$$
is a chain from $a$ to $d$ which represents $Ab$. Hence $a\not\preccurlyeq Ab$
by Lemma \lemChain(a).
\qed\enddemo


The rest of this section is devoted to the proof of Theorem \thAQP(b).
So, we fix $x,y\in\Cal A$ and $k,l\ge1$ and we
define $\Cal Q_m$,  $\Qlen(X)$, and $\Cal Q_{\min}(X)$ by (\eqQm)--(\eqQmin), see \S\sectProofOne.
We set also
$$
    \Cal Q_m^0 = \{ X\in\Cal Q_m\mid \ell(X)\le2m+k+l-2\}
$$
and $\Qmin^0(X)=\Cal Q^0_n\cap X^G$ for $n=\Qlen(X)$.
If $\Qlen(X)=0$, then the conclusion of Theorem 2(b) holds by definition of 
$\Qlen(X)$, so we shall consider the case when $\Qlen(X)>0$.

From now on, $x_1,x_2,\dots$ and $y_1,y_2,\dots$ will always denote some atoms which
are conjugate to $x$ and $y$ respectively.

\proclaim{ Lemma \lemAQ } {\rm(a)}. If $X\in\Cal Q_m$ and $m>0$,
then $\ell(X)\le 2m+k+l-1$.

\smallskip\noindent
{\rm(b)}. If $\Cal Q_m\cap X^G\ne\varnothing$ and $m>0$,
then $\Cal Q_m^0\cap X^G\ne\varnothing$.
In particular, $\Qmin^0(X)\ne\varnothing$ when $\Qlen(X)>0$.
\endproclaim

\demo{ Proof } (a). Let $X=P^{-1}x_1^k Py_1^l$, $\ell(P)=m$. We have $\ell(y_1^l)=l$
and, by Lemma \lemAQPone, we have $\ell(P^{-1}x_1^k P)\le 2m+k-1$.
Thus the result follows from
Lemma \lemTrIneq.

\smallskip
(b). Let $X_0=P^{-1}x_1^k P$ and $X=X_0 y_1^l$
with $\inf P=0$, $\ell(P)\le m$. We have to show that
$\Cal Q^0_m\cap X^G\ne\varnothing$.
By Lemma \lemAQPone, we may assume that the left normal form of $X_0$ is as
stated in Theorem \thAQP(a) with $n\le m$. If $n<m$, then the result follows from (a).
So, we suppose that $n=m$. Without loss of generality we may assume that
$P=B_1\dots B_n$. If $y_1\not\in F(\varphi(X_0))$, then $\varphi(X_0)y_1$
is a simple element by Lemma \lemBrSa, hence
$\ell(X)\le\ell(X_0)+\ell(y_1^l)-1$ and we are done. 
So, we suppose that $y_1\in F(\varphi(X_0))$.

\smallskip
Case 1. $m\ge2$. We have $\varphi(X_0)=B_n=B'_n y_1$, $B'_n\in\Cal P$.
Since $B_n$ is square free, we have $B'_n\not\succcurlyeq y_1$.
Let $P'=B_1\dots B_{n-1}B'_n$, $X'_0=(P')^{-1}x_1 P'$, $X'=X'_0y_1^l$.
Then we have $\ell(P')\le m$, hence $X'\in\Cal Q_m$.
The condition (\eqArtinAB) for $i=n$ can be rewritten as
$\Delta=B_n\tau^n(A_n)=B'_ny_1\tau^n(A_n)$, thus $A'_n=\tau^{-n}(y_1)A_n$
is a simple element such that $A'_n\Delta^{n-1}B'_n=\Delta^n$. Hence
$X'_0=\Delta^{-n}\cdot A'_n\cdot A_{n-1}\cdot \dots\cdot A_1\cdot x_1^{k-1}
\cdot x_1B_1\cdot B_2\cdot\dots\cdot B_{n-1}\cdot B'_n$ and we obtain
$X'=\Delta^{-n}\cdot A'_n\cdot A_{n-1}\cdot \dots\cdot A_1\cdot x_1^{k-1}
\cdot x_1B_1\cdot B_2\cdot\dots\cdot B_{n-1}\cdot B'_ny_1\cdot y^{l-1}$.
The number of simple factors in this decomposition is $2m+k+l-2$.
Thus $X'\in\Cal Q_m^0$. It remains to note that
$X' = y_1 X y_1^{-1}\in X^G$.

\smallskip
Case 2. $m=1$. We have $\varphi(X_0)= x_1B_1\succcurlyeq y_1$.
If $B_1\succcurlyeq y_1$, then we repeat the same arguments as in Case 1.
If $B_1\not\succcurlyeq y_1$, then the ``right-to-left version''
of Lemma \lemAG\ implies $B_1y_1=y_2B_1$, hence 
$X=B_1^{-1} x_1^k B_1 y_1^l = B_1^{-1} x_1^k y_2^l B_1\in\Cal Q_0^G$
which contradicts the
condition $\Qlen(X)=1$.
\qed\enddemo

\proclaim{ Lemma \lemAQi }
If $X\in\Qmin^0(X)$ and $\Qlen(X)>0$,
then the left normal form of $X$ is as stated in Theorem {\rm\thAQP(b)}
with $n=\Qlen(X)$.
\endproclaim

\demo{Proof}
Let $X\in\Qmin^0(X)$. Then $X=P^{-1}x_1^k Py_1^l$ with $\ell(P)=n=\Qlen(X)$.
Without loss of generality we may assume that $\inf P=0$ and
$\|P\|$ is the minimal possible among all
presentations of $X$ in this form.
Let $P=B_1\cdot\dots\cdot B_n$ be the left normal form of $P$ and
let $A_1,\dots,A_n$ be defined by (\eqArtinAB).
Then (\eqASC) represents $X$.

Case 1. $n\ge2$.
Let us show that (\eqASC) is left weighted and (\eqArtinAx), (\eqArtinBx) hold.
By Lemma \lemAQPone, the part $\Delta^{-n}\cdot A_n\cdot\dots\cdot B_n$ of
(\eqASC) is left weighted and (\eqArtinAx) holds 
(here we use the minimality of $\|P\|$).
So, it remains to prove that:
(i) $B_n y_1$ is a simple element;
(ii) (\eqArtinBx) holds;
(iii) $B_n y_1\cdot y_1$ is left weighted;
(iv) $B_{n-1}\cdot B_n y_1$ is left weighted.
Indeed:
\smallskip

(i). Otherwise $A_n\cdot\dots\cdot B_n\cdot y_1^l$ is left weighted, hence
$\ell(X)=2n+k+l-1$ which contradicts the fact that $X\in\Qmin^0(X)$.

\smallskip
(ii). Combine (i), Lemma \lemSBeqRA, and the fact that $B_n\tau^n(A_n)=\Delta$.

\smallskip
(iii). Follows from Lemma \lemBrSaSqF.

\smallskip
(iv).  Suppose that there exists $z\in R(B_{n-1})\cap S(B_n y_1)$.
Since $B_{n-1}\cdot B_n$ is left weighted, we have $z\not\in S(B_n)$.
Hence, by Lemma \lemAG, we have $B_n y_1^l = z^l B_n$. Thus $z\sim y_1$
and $B_n X {B_n}^{-1} = Q^{-1}x_1^k Q z^l$ where $Q=B_1\dots B_{n-1}$.
Since $\ell(Q)=n-1$, this contradicts the fact that $n=\Qlen(X)$.

\smallskip

Case 2. $n=1$. In this case (\eqASC) takes the form
$\Delta^{-1}\cdot A_1\cdot x_1^{k-1}\cdot x_1B_1y_1\cdot y_1^{l-1}$.
We have to show that this product is left weighted and (\eqArtinABx) holds,
that is:
(i) $x_1 B_1 y_1$ is a simple element;
(ii) (\eqArtinABx) holds;
(iii) $x_1 B_1 y_1\cdot y_1$ is left weighted;
(iv) $\varphi(A_1 x_1^{k-1})\cdot x_1 B_1 y_1$ is left weighted;
(v) $A_1\cdot x_1$ is left weighted.
Indeed:

\smallskip
(i). Otherwise $A_1\cdot x_1^{k-1}\cdot x_1B_1\cdot y^l$ is
left weighted (because $A_1\cdot x_1^{k-1}\cdot x_1B_1$ is so by Lemma \lemAQPone),
hence $\ell(X)=k+l+1$ which contradicts the fact that $X\in\Qmin^0(X)$.

\smallskip
(ii). Let $\tilde y_1=\delta^{-1}(y_1)$ (as in (\eqArtinABx)).
By (i) we have $y_1\in R(B_1)=S(\tau(A_1))$ and $x_1\in L(B_1)=F(A_1)$.
So, $\tilde y_1\preccurlyeq A_1=A'_1 x_1$ with $A'_1\in\Cal P$.
We have to show that $\tilde y_1\preccurlyeq A'_1$.
Suppose that $\tilde y_1\not\preccurlyeq A'_1$.
Then it follows from Lemma \lemAG\ that $\tilde y_1 A'_1=A'_1 x_1$.
Hence $y_1\sim x_1$ and $\tilde y_1 A_1 = \tilde y_1 A'_1 x_1 = A'_1 x_1^2 = A_1 x_1$.
Thus $X=\Delta^{-1} A_1 x^k B_1 y_1^l =\Delta^{-1} \tilde y_1^k A_1 B_1 y_l = y_1^{k+l}\in\Cal Q_0$
which contradicts the fact that $X\in\Qmin^0(X)$.

\smallskip
(iii). Follows from Lemma \lemBrSaSqF.

\smallskip
(iv). Suppose that there exists $z\in R(\varphi(A_1x_1^{k-1}))\cap S(x_1 B_1 y_1)$.
Since $\varphi(A_1x_1^{k-1})\cdot x_1B_1$ is left weighted by Lemma \lemAQPone, we have
$z\preccurlyeq x_1 B_1 y_1$ and $z\not\preccurlyeq x_1 B_1$. By Lemma \lemAG,
it follows that $z x_1 B_1 = x_1 B_1 y_1$. Hence, $z\sim y_1$ and
$X=B^{-1} x_1^{k-1} (x_1 B_1) y_1^l = B_1^{-1} x_1^{k-1} z^l (x_1 B_1) \sim x_1^k z^l\in\Cal Q_0$
which contradicts the fact that $X\in\Qmin^0(X)$.

\smallskip
(v). Combine (\eqArtinABx) and Lemma \lemBrSaSqF.
\qed\enddemo

\proclaim{ Lemma \lemAQii }
If $X\in\Qmin^0(X)$ and $\Qlen(X)>0$, then $\frak s(X)\in\Qmin^0(X)$.
\endproclaim

\demo{Proof}
By Lemma \lemAQi, we may assume that the left normal form of $X$ is as stated in
Theorem 2(b) with $n=\Qlen(X)$.

\smallskip
Case 1. $n\ge 2$ or $l\ge 2$.
Let $\tilde A_n = \iota(X)=\tau^n(A_n)$
and $Y=\tilde A_n^{-1} y_1^l \tilde A_n = \Delta^{-1} B_n y_1^l\tilde A_n$.
By Lemma \lemAQPone, the left normal form of $Y$ is
$\Delta^{-1}\cdot B'_n\cdot y_2^{l-1}\cdot y_2\tilde A'_n$ where $B'_n$ and $\tilde A'_n$
are simple elements such that $B'_n\tilde A'_n=\Delta$ and $B'_n=B''_n y_2$, 
$B''_n\in\Cal P$.
We can rewrite the left normal form of $Y$ also as
$\Delta^{-1}\cdot B''_n y_2\cdot y_2^{l-1}\cdot\tilde A''_n$ where
$\tilde A''_n=y_2\tilde A'_n$. Let $A''_n=\tau^{-n}(\tilde A''_n)$.
Then, by Corollary \corCharneyR(a) (if $n>1$) or by Corollary \corCharneyR(b) (if $n=1$
and $l>1$), we have
$$
   \frak s(X) = \Delta^{-n}\cdot A''_n\cdot A_{n-1}\cdot\dots\cdot A_1\cdot x_1^k\cdot
              B_1\cdot B_2\cdot\dots\cdot B_{n-1}\cdot B''_n \cdot y_2^l.
$$
Hence $\frak s(X)\in\Cal Q_n$. Since, $\ell(\frak s(X))\le\ell(X)$
(see [\refGebGM; Lemma 1]), we conclude that $\frak s(X)\in\Qmin^0(X)$.

\smallskip
Case 2. $n=l=1$. Combining (\eqASC) and (\eqArtinABx) and denoting $A''_1$ by $A$
and $B_1$ by $B$, we may rewrite the left
normal form of $X$ a more symmetric way as
$\Delta^{-1}\cdot\tilde y_1A x_1\cdot x_1^{k-1} \cdot x_1 B y_1$ where
$\tilde y_1 A x_1 B = A x_1 B y_1 = \Delta$ (and hence $\tau(\tilde y_1)=y_1$).
Then we have
$\bold d(X) =  \tilde x_1\tilde B\tilde y_1 \cdot\tilde y_1 A x_1\cdot x_1^{k-1}$
where $\tau(\tilde x_1) = x_1$ and  $\tau(\tilde B) = B$.

Let us define
$\bar{\Cal Q}_m$, $\bar{\Cal Q}_m^0$, etc. in the same way as
$\Cal Q_m$, $\Cal Q_m^0$, etc. but with $x^k$ and $y^l$ exchanged.
Then we have $\bold d(X)\in\bar{\Cal Q}_1^0$.
It is clear that $\Cal Q_m\cap X^G\ne\varnothing$ if and only if
$\bar{\Cal Q}_m\cap X^G\ne\varnothing$. Since, moreover,
$\ell(\bold d(X))\le\ell(X)$, we conclude that
$\bold d(X)\in\bar{\Cal Q}^0_{\min}(X)=\bar{\Cal Q}_1^0\cap X^G$.
Then, by Lemma \lemAQi\ applied to $\bar{\Cal Q}_{\min}^0(X)$,
the left normal form of $\bold d(X)$
is $\Delta^{-1}\cdot\tilde x_2\tilde B'\tilde y_2\cdot \tilde y_2 A'
x_2\cdot x_2^{k-1}$ where
$\tilde x_2\tilde B'\tilde y_2 A' = \tilde B'\tilde y_2 A' x_2 = \Delta$. Hence
$\bold c(\bold d(X)) = \Delta^{-1}\cdot\tilde y_2 A' x_2\cdot x_2^{k-1}\cdot x_2 B' y_2
\in\Qmin^0(X)$.
where $B'=\tau(\tilde B')$ and $y_2=\tau(\tilde y_2)$. It remains to note
that $\bold c(\bold d(X)) = \frak s(X)$ by Lemma \lemCDcommut.
\qed\enddemo

Theorem 1(b) follows immediately from Lemma \lemAQ(b), Lemma \lemAQi,
and Lemma \lemAQii\ combined with the fact that $\frak s^m(X)\in\SC(X)$ for $m$
sufficiently large.

%
%

\head\sectExamples. An example
\endhead

It is shown in [\refOrevkovQPthree] that if a braid $X$ with three strings 
is quasipositive, then any positive word $W$ in the standard
generators $\sigma_1$, $\sigma_2$ of $\Br_3$ such that $X=\Delta^p W$ with $p\le 0$,
satisfies the following property. There exists a word $W'$ obtained
by removing $e(X)$ letters from $W$ such that $\Delta^p W'=1$.
The same result is true for the dual Garside structure on $\Br_3$.

Theorems 1 and 2 of the present paper show that if $X$ is a quasipositive braid
with any number of strings but with $e(X)\le 2$, then $\SC(X)$ contains an
element which can be presented in the form $\Delta^p W$ where $W$ is a positive
word which satisfies the above property.

The following example shows that this is no longer true in the dual Garside structure on
$\Br_4$ for braids of algebraic length 3. Namely, let $\sigma_1$, $\sigma_2$, $\sigma_3$
still denote the standard Artin generators of $\Br_4$.
Let $\delta=\sigma_3\sigma_2\sigma_1$,
$\sigma_0=\sigma_3^\delta$, $\alpha=\sigma_1^{\sigma_2}$, $\beta=\sigma_2^{\sigma_3}$.
Then $\sigma_0,\dots,\sigma_3,\alpha,\beta$ are the atoms and $\delta$
is the Garside element of the Birman-Ko-Lee Garside
structure [\refBKLone] on $\Br_4$. Let
$$
  X = \delta^{-1}\cdot \beta\cdot\alpha\cdot\sigma_1
       \cdot\sigma_2\cdot\alpha\cdot\beta
                                                                        \eqno(\eqExOne)
$$
This braid is quasipositive, indeed, if we remove the second $\alpha$,
then we obtain
$$
      \delta^{-1}\cdot \beta\cdot\alpha\cdot\sigma_1\cdot\sigma_2\beta
    = \delta^{-1}\cdot \beta\cdot\alpha\cdot\sigma_1\cdot\sigma_3\sigma_2
    = \delta^{-1}\cdot \beta\cdot\alpha\cdot\sigma_1\sigma_3\cdot\sigma_2
$$
which is of the form (\eqSC) with $n=1$, $x_1=\alpha$, $y_1=\sigma_2$, $A_1=\beta$,
$B_1=\sigma_1\sigma_3$.
The braid $X$ is rigid and (\eqExOne) is its left (and also right) normal form,
so, $\bold c^6(X)=\tau(X)$. The cycling orbit of $X$ contains 24 elements and it can
be easily checked that it coincides with the summit set $\SS(X)$ (and hence with
$\SSS(X)$, $\USS(X)$, and $\SC(X)$).
Thus, for any presentation of any element of $\SS(X)$ in the form $\delta^{-1}W$
with a positive word $W$, it is impossible to remove three letters from
$W$ to obtain the trivial braid.

%
%

\head \sectBthree. Quasipositivity problem for 3-braids
\endhead

The result of [\refOrevkovQPthree] cited in \S\sectExample\
leads to an evident algorithm to decide if a given 3-braid $X$ is
quasipositive or not: it is enough to try to remove $e(X)$ letters
from $W$ in all possible ways.
Here we give a minor improvement of this algorithm in the 'branch and bound' style.
The new algorithm is still of exponential time with respect to the algebraic length $e(X)$
but the base of the exponent is smaller.
The improvements are based on the simple observations
summarized in Proposition \propTi\ below.

Given $\vec a=(a_1,\dots,a_n)$, $a_i>0$, and $p\in\Bbb Z$, we set
$\len(\vec a)=n$ and
$$
    X(p,\vec a)=X(p;a_1,\dots,a_n) = \Delta^p\,
     \underset\text{$n$ alternating factors}\to{
      \underbrace{
     \sigma_1^{a_1}\sigma_2^{a_2}\sigma_1^{a_3}\sigma_2^{a_4}\dots}}\;\in\Br_3  \eqno(\eqTi)
$$
We say that $(p',\vec a')$,
is obtained from $(p,\vec a)$ by an {\it elementary reduction}
in the following cases:
\roster
\item"(R1)" $n\ge 2$, $n\not\equiv p\mod 2$, $p'=p$, $\vec a'=(a_1+a_n,a_2,\dots,a_{n-1})$;
\item"(R2)" $n\ge 3$, $a_2=1$, $a_1,a_3\ge2$, $p'=p+1$, $\vec a'=(a_1-1,a_3-1,a_4,\dots,a_n)$;
\item"(R3)" $p$ is even, $\vec a=(1,a_2)$, $a_2\ge3$, $p'=p+1$, $\vec a'=(a_2-2)$
\item"(R4)" $p$ is even, $\vec a=(1,2)$, $p'=p+1$, $\vec a'=(\;)$;
\item"(R5)" $n\ge4$, $a_2=a_3=1$, $p'=p+1$, $\vec a'=(a_1+a_4-1,a_5,\dots,a_n)$;
\item"(R6)" $p$ is odd, $\vec a=(1,1,a_3)$, $a_3\ge2$, $p'=p+1$, $\vec a'_1=(a_3-1)$;
\item"(R7)" $p$ is odd, $\vec a=(1,1,1)$, $p'=p+1$, $\vec a'=(\;)$;
\item"(R8)" $p\equiv n\mod2$ and $(p',\vec a')$ is obtained from $(p,\vec a)$ by
a cyclic permutation of $\vec a$ followed by one of (R2)--(R7).
\endroster
A pair $(p,\vec a)$, $n=\len(\vec a)$, is called {\it reduced} if no elementary reduction 
can be applied to it. This is equivalent to the fact that either
$$
   \text{(i) }n\le1,\quad
   \text{or (ii) }\vec a=(1,1),\; p\equiv0(2),\quad\text{or (iii) }
\text{$n\equiv p(2)$ and all $a_i\ge2$.}             \eqno(\eqTii)
$$
It is clear that if $(p',\vec a')$ is an elementary reduction of $(p,\vec a)$, then
$X(p',\vec a')$ is conjugate to $X(p,\vec a)$. It follows easily from the Garside theory that
if a pair $(p,\vec a)$ is reduced, then $\inf_s X(p,\vec a)=p$
and $(p,\vec a)$ is determined by the conjugacy class of $X(p,\vec a)$
up to cyclic permutation of $\vec a$.

\proclaim{ Lemma \lemT }
Suppose that $(p',\vec a')$ is obtained from $(p,\vec a)$ by an elementary reduction.
Let $n=\len(\vec a)$, $n'=\len(\vec a')$.
Then $p'+n' \le p+n$. \qed
\endproclaim

For a braid $X$, we denote the signature and the nullity of its closure by $\Sign(X)$ and $\Null(X)$
respectively (we follow the convention that the nullity of a link is the nullity
of the symmetrized Seifert form corresponding to a
connected Seifert surface).
If a braid $X$ is quasipositive, then Murasugi-Tristram inequality implies
$$
    1+\Null(X) \ge |\Sign(X)|+m-e(X)                   \eqno(\eqMT)
$$
where $m$ is the number of strings (see details in [\refOrevkovGAFA; \S3.1]).
The following fact can be easily derived from [\refOrevkovDN; Prop.~8.2] or from [\refGG; Th.~4.2]
(also it was conjectured and partially proved in [\refMurasugi; \S\S9--11]).

\proclaim{ Lemma \lemTsign } Let $X=X(p,\vec a)$ with $(p,\vec a)$ reduced, $n=\len(\vec a)\ge 2$,
and $\vec a\ne(1,1)$. Then $\Sign(X)+\Null(X)=p+n-e(X)$ and
$$
    \Null(X) = \cases 1, &\text{if }\vec a=(2,2,\dots,2)
    \text{ and }p+n\equiv0\mod4,\\ 0,&\text{otherwise.}           \qed
    \endcases
$$
\endproclaim

Let us denote a sequence $(2,2,\dots,2)$ ($n$ times) by $2_n$.

\proclaim{ Lemma \lemTi }
{\rm(a)}. If $q$ is even and $n\ge0$, then $\Delta^q\sigma_1^{-n}\sim X(q-n;2_n)$.
If $q$ is odd, then $\Delta^q\sigma_1^{-1}\sim X(q-1;1,1)$,  $\Delta^q\sigma_1^{-2}\sim X(q-1;1)$,
and $\Delta^q\sigma_1^{-k}\sim X(q-k+1;3,2_{k-3})$ for $k\ge 3$.

\smallskip
{\rm(b).} A braid $X=\Delta^q\sigma_1^{-n}$ is quasipositive if and only if
either $(q,n)=(0,0)$, or $q\ge0$ and $2n<5q$.
\endproclaim

\demo{ Proof }
(a). Evident.

\smallskip
(b). Let $q$ be even and $n\ge 2$. Then $X$ is quasipositive
if and only if $X'=\Delta^{q-1}\sigma_1^{-(n-2)}$ is quasipositive.
Indeed, by (a), we have $X\sim X(q-n;2_n)$, hence, by Proposition \propTiReduc, $X$ is
quasipositive if and only if one of $X_i=X(q-n,f_i(2_n))$ is. For any $i$ we have
$X_i\sim X(q-n;2_2,1,2_{n-3})
    \overset\text{(R2)}\to\sim X(q-n+1;2,1,1,2_{n-4})
    \overset\text{(R5)}\to\sim X(q-n+2;3,2_{n-5}) \sim X'$ (we suppose here that $n\ge5$
and we leave to the reader to check that $X_i\sim X'$ for $n=2,3,4$).
Since $q$ is even, we have $2n<5q\Leftrightarrow2n<5q-1\Leftrightarrow2(n-2)<5(q-1)$,
thus it is enough to prove the statement only for odd $q$. From now on we suppose that $q$ is odd.

\smallskip
Suppose that $0<2n<5q$.
Let us prove by induction that $X$ is quasipositive. 
If $q=1$, then $n\le2$
and $X=\Delta\sigma_1^{-2}=\sigma_1\sigma_2\sigma_1^{-1}$ is quasipositive.
If $q\ge 3$, then we have $\Delta^q\sigma_1^{-2}\sigma_2^{-1}\sigma_1^{2-n}=
\Delta^q\sigma_1^{-1}\Delta^{-1}\sigma_1^{3-n}=\Delta^{q-1}\sigma_2^{-1}\sigma_1^{3-n}
=\sigma_1\Delta^{q-1}\sigma_1^{-1}\sigma_2^{-1}\sigma_1^{3-n}
=\sigma_1(\Delta^{q-2}\sigma_1^{5-n})\sigma_1^{-1}$, hence
$$
  \Delta^q\sigma_1^{-n}=(\sigma_2^{-2}\sigma_1\sigma_2^2)
   \cdot\Delta^q\sigma_1^{-2}\sigma_2^{-1}\sigma_1^{2-n}=
  (\sigma_2^{-2}\sigma_1\sigma_2^2)
   \cdot\sigma_1(\Delta^{q-2}\sigma_1^{5-n})\sigma_1^{-1}.
$$
So, if
$\Delta^{q-2}\sigma_1^{5-n}$ is quasipositive by the induction hypothesis, then $X$ is also.

\smallskip
Suppose that $X$ is quasipositive.
Then (\eqMT) combined with (a) and with Lemma \lemTsign\ yields $2n\le 5q-1$.
\qed\enddemo

\definition{ Remark \remT }
In [\refOrW], the question of the quasipositivity of $X=\Delta^q\sigma_1^{-n}\in\Br_k$
is studied for any $k$. In particular, it is shown that this is so for
$n\le qk^2/3+O(qk)$. 
However, for $k=3$, the construction from [\refOrW] gives the quasipositivity
of $X$ only when $n\le 2q$ which is weaker than Lemma \lemTi(b).
\enddefinition

Given $\vec a=(a_1,\dots,a_n)$ and $i\in\{1,\dots,n\}$, we set
$f_i(\vec a)=(a'_1,\dots,a'_n)$ where $a'_i=1$ and $a'_j=a_j$ for $j\ne i$.

\proclaim{ Proposition \propTi } Let $(p,\vec a)$, $n=\len(\vec a)$, satisfy {\rm(\eqTii)}.
Let $X=X(p,\vec a)$. Then:

\smallskip\noindent
{\rm(a)}. If $p\ge 0$, then $X$ is quasipositive.

\smallskip\noindent
{\rm(b)}. If $p<0$ and $X$ is quasipositive, then 
$$
                 0<p+n<2e(X).                      \eqno(\eqTiii)
$$

\smallskip\noindent
{\rm(c)}. If $3n+5p>0$, then $X$ is quasipositive {\rm(}see Figure 1\/{\rm)}.

\smallskip\noindent
{\rm(d)}. If $3n+5p=0$ and $\vec a\ne(2_n)$, then $X$ is quasipositive.
Note that {\rm(\eqTiii)} implies $\vec a\ne(2_n)$ when  $3n+5p=0$,

\smallskip\noindent
{\rm(e)}. $X$ is quasipositive if and only if there exists $i$ such that the braid
$X(p,f_i(\vec a))$ is quasipositive.

\smallskip\noindent
{\rm(f)}. Suppose that $X(p,f_i(\vec a))$ is not quasipositive and $(p',\vec a')$ is obtained
from $(p,\vec a)$ by an elementary reduction. If $a_i$ is not involved in the reduction and
$a'_{i'}$ is the entry of $\vec a'$ which corresponds to $a_i$, then $X(p',f_{i'}(\vec a'))$
is not quasipositive.
\endproclaim

\midinsert

\centerline{\epsfxsize=55mm\epsfbox{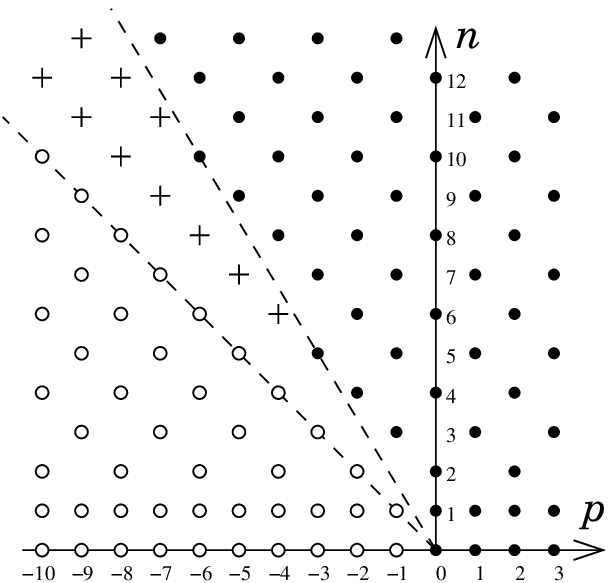}\hskip-5mm
\hbox to 65mm{\vbox{\hsize=65mm
\noindent
\roster
\item"$\bullet$" means that $X(p,\vec a)$ with given $p$ and $n=\len(\vec a)$ is quasipositive
if $(p,\vec a)$ is reduced; when $3n+5p=0$, we assume also that $\vec a\ne(2_n)$
\medskip
\item"$\circ$"  means that any $X(p,\vec a)$
with given $p$ and $n=\len(\vec a)$ is not quasipositive if $(p,\vec a)$ is reduced
\medskip
\item"$+$" $X(p,\vec a)$ may or may not be quasipositive
\endroster
\bigskip\medskip
}}}
\botcaption{ Figure 1 } When $X(p,\vec a)$ is {\sl a priori} (non-)quasipositive \endcaption
\endinsert

\demo{ Proof }
(a) and (f). Evident.

\smallskip
(e) Suppose that $X$ is quasipositive. If $p>0$, then the statement is obvious.
Suppose that $p\le0$. Then, by [\refOrevkovQPthree; Prop.~3.1],
one can remove some letters from the positive part of the right hand side of (\eqTi) so that
the resulting braid becomes trivial. This means that there exists $\vec a'=(a'_1,\dots,a'_n)$
such that $a'_i\le a_i$ for each $i$ and $X(p,\vec a')=1$.
It remains to note that if $a'_i\ge 2$ for all $i$, then $X(p,\vec a')\ne1$.

\smallskip
(b). Suppose that $p<0$ and $X$ is quasipositive. Then $n\ge 2$ by (e).
Thus, when $\vec a\ne(2_n)$, the result follows from (\eqMT) combined with Lemma \lemTsign.
If  $\vec a=(2_n)$, then the result follows from Lemma \lemTi.
Note that the left inequality $0<p+n$ can be proven also by induction using (e) and Lemma \lemT.

\smallskip
(c). It is clear that if $\vec a'=(a'_1,\dots, a'_n)$ is such that $a'_i\le a_i$ for any $i$, then
the quasipositivity of $X(p,\vec a')$ implies that of $X(p,\vec a)$. Thus, the result follows
from Lemma \lemTi\ if we set $\vec a'=(2_n)$.

\smallskip
(d). The same proof but with $\vec a'=(3,2_{n-1})$.
\qed\enddemo

Thus we obtain the following recursive algorithm.
The input is the pair $(p,\vec a)$ together
with the information about the indices $i$ for which it is known already that $X(p,f_i(\vec a))$
is not quasipositive, see Proposition \propTiBB.
The pair $(p,\vec a)$ is assumed to be {\it almost reduced}, i.~e., $p\equiv n=\len(\vec a)\mod 2$
when $n\ge2$, $a_0\ge1$, and $a_i\ge2$ for $i>0$ (since the algorithm is implemented below in {\tt C}
programming language, we assume here
that the entries of $\vec a$ are numbered from $0$ to $n-1$). First we
reduce $(p,\vec a)$ 
and check if the conclusion can be done using Proposition \propTi(a--d).
Then we check recursively if any of $X(p,f_i(\vec a))$ is quasipositive (see Proposition \propTiReduc)
taking into account the information that some of them are already known not to be.

Below we present an implementation of this algorithm in the form of
a {\tt C} function {\bf qp3()}. We assume that the input braid
is given in the form (\eqTi) with $(p,\vec a)$ almost reduced (the arguments {\bf p} and {\bf a}).
The argument {\bf n} should be equal to $\len(\vec a)$ and the argument {\bf e}
should be equal to $e(X(p,\vec a))=3p+\sum a_i$.
We assume that the pointer {\bf a} points to a preallocated array of at least
{\tt 2*e*n} integers. The first $n$ entries of this array contain the vector $\vec a$
and the others are used for the intermediate data. The initial values of the array
will be lost after the computation.

During the computations, we assume that the
vector $\vec a$ is represented by the absolute values of the entries of the array {\bf a}
whereas the negative sign of {\bf a\tt[\bf i\tt]} is used to encode the information that the
braid $X(p,f_i(\vec a))$ is not quasipositive, see Proposition \propTiBB.
Instead of computing $f_i(\vec a)$ for $i\ge1$, we compute $f_0$ of cyclic permutations of $\vec a$
(this ensures that the input is always almost reduced).
The function {\bf qp3()} returns 1 if $X(p,\vec a)$ is quasipositive and 0 otherwise.
\medskip
\chardef\lb=`\{
\chardef\rb=`\}
\chardef\am=`\&
\chardef\us=`\_
{\tt
int qp3( int p, int *a, int n, int e )\lb\par
\  \  while( n>1 )\lb\ \ \ \ \ // reduce (p,a) assuming abs(a[i])>1 for i>0\par
\  \  \  \  if( a[0]==1 \am\am\ a[n-1]==1 )\lb\par
\  \  \  \  \  \  if( n==2 )break; else p++;\par
\  \  \  \  \  \  if( n==3 )\lb\ a[0]=abs(a[1])-1; n=1; break; \rb\par
\  \  \  \  \  \  a[1]=abs(a[1])+abs(a[n-2])-1; a[0]=a[n-=3]; break; \rb\par
\  \  \  \  if( a[0]==1 )a++; else\lb\ if( a[n-1]!=1 )break; \rb\par
\  \  \  \  p++;n--;a[0]=abs(a[0])-1;a[n-1]=abs(a[n-1])-1;\rb\ \ \ // reduced\par
\  \  if( p >= 0 )return 1;\ \ \ \ \ \ \ \ \ \ \ \ \ \ \ \ \ \ \ \ // see Prop.~\propTi(a)\par
\  \  if( !(0 < p+n \am\am\ p+n < 2*e) )return 0;\ \ \             // see Prop.~\propTi(b)\par
\  \  if( 3*n + 5*p >= 0 )return 1;\ \ \ \ \ \ \ \ \ \ \ \         // see Prop.~\propTi(c,d)\par
\  \  \lb\  int count=n,e1,*a1,i;\par
\  \  \  \  while( count-- )\lb\ // repeat n times\par
\  \  \  \  \  \  if( a[0] > 0 )\lb\ // a[0]<0 means that X(p,f\us0(a)) is not qp\par
\  \  \  \  \  \  \  \  if( (e1=e-a[0]+1) >= 0 )\lb\par
\  \  \  \  \  \  \  \  \  \  for( a1=a+n,i=1; i<n; i++ )a1[i]=a[i];\par
\  \  \  \  \  \  \  \  \  \  a1[0]=1; if( qp3(p,a1,n,e1) )return 1; \rb\ // recursion\par
\  \  \  \  \  \  \  \  a[0]=-a[0]; \rb\par
\  \  \  \  \  \  a[n] = (*a++); \rb\ \ \ \ \ // cyclic permutation of the array a\par
\  \  \  \  return 0; \rb\rb
}

%
%

\head \sectBlock. Blocking property of the dual Garside structures on Artin groups
\endhead

In this section we prove a property (we call it the {\it blocking property}\/)
of square free symmetric homogeneous Garside structures, in particular, the dual Garside
structures on Artin groups and the Garside structure [\refBC] on $G(e,e,r)$.
This property is not used in this paper but we hope it to be useful for
the quasipositivity problem in the general case.

\proclaim{ Proposition \propBlock } Let $(G,\Cal P,\delta)$ be a square free
symmetric homogeneous Garside structure. Let
$k\ge1$, $A\in{]1,\Delta[}$, $B=\partial A$ {\rm(}i.~e., $AB=\delta${\rm)}
and let $x$ be an atom such that
$X=A\cdot x^k\cdot B$ is in left normal form.
Let $Y\in G$, $\inf Y=0$.
Then either $\delta\preccurlyeq XY$, or $\iota(XY)=A$.
\endproclaim

\proclaim{ Corollary \corBlock } Let $(G,\Cal P,\delta)$ be a square free
symmetric homogeneous Garside structure and let $X$ be as in Theorem 1(a), thus
the left normal form of $X$ is given by {\rm(\eqQPone)} and {\rm(\eqAB)}.
Let $Y\in G$, $\inf Y=0$. Then either $\inf(XY)>\inf X$, or the left normal form of $XY$
begins with $\delta^{-n}\cdot A_n\cdot\dots\cdot A_1$.
\endproclaim

\proclaim{ Lemma \lemBi } {\rm(Compare with [\refBKLone; Cor.~3.7])}.
Let $(G,\Cal P,\delta)$ be a square free and symmetric
Garside structure. Let $A$ be a simple element of $G$
and let $S(A)=\{x_1,\dots,x_m\}$. Then
$A = x_1\vee\dots\vee x_m$.
\endproclaim

\demo{ Proof } Let $B = x_1\vee\dots\vee x_m$. Then $B\preccurlyeq A$, i.~e.,
$A=BC$ for $C\in[1,\Delta]$. We have to prove that $C=1$. Suppose that $C\ne1$.
Let $y\in S(C)$. Since $A\succcurlyeq C$ and the Garside structure is symmetric,
we have $C\preccurlyeq A$, hence $y\preccurlyeq C\preccurlyeq A$, i.~e., $y\in S(A)$.
Hence $y\preccurlyeq B$ by the definition of $B$. Since the Garside structure is symmetric,
it follows that $B\succcurlyeq y$. Thus we have $y\in F(B)$ and $y\in S(C)$ which contradicts the
fact that $A=BC$ is square free.
\qed\enddemo

\proclaim{ Lemma \lemBii } Let $(G,\Cal P,\delta)$ be a
symmetric homogeneous Garside structure. Let $x$ and $y$ be atoms such that
$xy\not\preccurlyeq\delta$. Let $D=x^{-1}(x\vee y)$, Then $y\vee D=x\vee y$.
\endproclaim

\demo{ Proof }
We have $x\vee y=xD$. Since the Garside structure is symmetric, it follows that
$D\preccurlyeq x\vee y$, hence $y\vee D\preccurlyeq x\vee y$.
Since the Garside structure is homogeneous,
it follows that $\|x\vee y\|=\|xD\|=\|D\|+1$ and we obtain
$$
      D\preccurlyeq y\vee D \preccurlyeq x\vee y
\quad\text{and}\quad \|x\vee y\|=\|D\|+1.
$$
Thus, it is enough to show that $D\ne y\vee D$. Suppose that
$D=y\vee D$. Then we have $y\preccurlyeq D$, hence
$xy\preccurlyeq xD=x\vee y\preccurlyeq\delta$.
Contradiction.
\qed\enddemo

\proclaim{ Lemma \lemBiii }
Let $(G,\Cal P,\delta)$ be a
symmetric square free Garside structure.
Let $A\in[1,\Delta]$ and $P\in\Cal P$.
Then $\iota(A^2P)= \iota(AP)$. In particular, $S(A^2P)=S(AP)$.
\endproclaim

\demo{ Proof } Let $B=\iota(AP)$. By Lemma \lemCharneyA, we have
$\iota(A^2P)=\iota(AB)$. We have $B=AC$ for a simple element $C$.
Since $B$ is simple and the Garside structure is symmetric, we have
$B=AC=CA'$ with $A'\in\Cal P$. Hence $AB=ACA'=BA'$.
We have $F(A')=S(A')$ (because the Garside structure is symmetric)
and $R(A')\subset\Cal A\setminus F(A')$ (because the Garside structure is square free;
$\Cal A$ stands for the set of atoms).
Hence $R(B)=R(CA')\subset R(A')\subset\Cal A\setminus F(A')=\Cal A\setminus S(A')$
which means that the decomposition $AB=B\cdot A'$ is left weighted.
Thus $\iota(A^2P)=\iota(AB)=\iota(B\cdot A')=B=\iota(AP)$
\qed\enddemo

\demo{ Proof of Proposition \propBlock }
Suppose that $A\ne\iota(Ax^kBY)$. Then $R(A)\cap S(x^kBY)\ne\varnothing$. Let
$y\in R(A)\cap S(x^kBY)$.
By Lemma \lemSBeqRA\ we have $R(A)=S(B)$, hence
$$
          y\in S(B).                                     \eqno(\eqBi)
$$
Let $D=x^{-1}(x\vee y)$.
Since $y\in S(x^kBY)$, we have $x\vee y\preccurlyeq x^kBY$, i.~e., $xD\preccurlyeq x^kBY$.
By Lemma \lemBiii, this implies  $xD\preccurlyeq xBY$.
By canceling $x$, we obtain $D\preccurlyeq BY$.  Combining this fact with (\eqBi),
we obtain
$$
   y\vee D\preccurlyeq BY.                                     \eqno(\eqBii)
$$
Combining (\eqBi) with the fact that
$A\cdot x^k\cdot B$ is left weighted,
we obtain $xy\not\preccurlyeq\delta$. Hence, by Lemma \lemBii, we have
$y\vee D=x\vee y$. Hence, by (\eqBii), we obtain
$$
   x\preccurlyeq x\vee y=y\vee D\preccurlyeq BY.                \eqno(\eqBiii)
$$

Let us prove that $B\preccurlyeq x^kBY$. By Lemma \lemBi, it is enough to show
that $S(B)\subset S(x^kBY)$. Let $z\in S(B)$ and let $E=x^{-1}(x\vee z)$.
Combining (\eqBiii) with $z\preccurlyeq B\preccurlyeq BY$, we obtain
$x\vee z\preccurlyeq BY$, i.~e., $xE=x\vee z\preccurlyeq BY$.
Since the Garside structure is symmetric and $xE\preccurlyeq\delta$,
it follows that $E\preccurlyeq xE\preccurlyeq BY$, hence, $xE\preccurlyeq xBY$
and we conclude that $z\preccurlyeq x\vee z=xE\preccurlyeq xBY$.
Thus we have proven that $S(B)\subset S(xBY)$. By Lemma \lemBiii, it follows that $S(xBY)=S(x^kBY)$,
hence $S(B)\subset S(x^kBY)$. By Lemma \lemBi, this implies
$B\preccurlyeq x^kBY$. Multiplying this inequality by $A$,
we obtain $\delta=AB\preccurlyeq Ax^kBY=XY$.
\qed\enddemo


\Refs
\def\r{\ref}

\r\no\refBessis
\by D.~Bessis
\paper The dual braid monoid
\jour Ann. Sci. \'Ecole Norm. Sup. \vol 36 \yr 2003 \issue 5 \pages 647--683
\endref

\r\no\refBessisCorran
\by D.~Bessis, R.~Corran
\paper Non-crossing partitions of type $(e,e,r)$
\endref

\r\no\refBGGMii
\by  J.~S.~Birman, V.~Gebhardt, J.~Gonz\'alez-Meneses
\paper Conjugacy in Garside groups II: structure of the ultra summit set
\jour  Groups, Geom. and Dynamics \vol 1 \yr 2008 \pages 13--61
\endref

\r\no\refBKLone
\by J.~Birman, K.-H.~Ko, S.-J.~Lee
\paper A new approach to the word and conjugacy problems in the braid groups
\jour  Adv. Math. \vol 139 \yr 1998 \pages 322--353
\endref

\r\no\refBKLtwo
\by J.~Birman, K.-H.~Ko, S.-J.~Lee
\paper The infimum, supremum, and geodesic length of a braid conjugacy class
\jour Adv. Math. \vol 164 \yr 2001 \pages 41--56
\endref

\r\no\refBrSa
\by E.~Brieskorn, K.~Saito
\paper Artin-Gruppen und Coxeter-Gruppen
\jour Invent. Math. \vol 17 \yr 1972 \pages 245--271
\endref

\r\no\refCharney
\by R.~Charney
\paper Artin groups of finite type are biautomatic
\jour Math. Ann. \vol 292 \yr 1992 \pages 671--683
\endref

\r\no\refCrisp
\by J.~Crisp
\paper Injective maps between Artin groups
\inbook Geometric Group Theory Down Under, Proceedings of a Special Year
in Geometric Group Theory (J.~Cossey et al, ed.)
\publ W. de Gruyter \yr 1999 \pages 119--137
\endref

\r\no\refDehornoy
\by P.~Dehornoy
\paper Groupes de Garside
\jour Ann. Sci. \'Ecole Norm. Sup. \vol 35 \yr 2002 \pages 267--306
\endref

\r\no\refDP
\by P.~Dehornoy, L.~Paris
\paper Gaussian groups and Garside groups, two generalizations of Artin Groups
\jour  Proc. London Math. Soc. (3) \vol 79 \yr 1999 \pages 569--604
\endref

\r\no\refDeligne
\by P.~Deligne
\paper Les immeubles des groupes de tresses g\'en\'eralis\'es
\jour Invent. Math. \vol 17 \yr 1972 \pages 273--302
\endref

\r\no\refEM
\by E.~ElRifai, H.~Morton
\paper Algorithms for positive braids
\jour  Quart. J. Math. Oxford Ser. (2) \vol 45 \yr 1994 \pages 479--497
\endref

\r\no\refThurstonAndCo
\by D.~B.~A.~Epstein, J.~W.~Cannon, D.~F.~Holt, S.~V.~F.~Levy, M.~S.~Paterson, W.~P.~Thurston
\book Word Processing in Groups
\publ Jones \& Bartlett \publaddr Boston, MA \yr 1992 \pages Chapter 9.
\endref

\r\no\refGG
\by J.-M.~Gambaudo, E.~Ghys
\paper Braids and signatures
\jour Bull. Soc. Math. France \vol 133 \yr 2005 \pages 541--579
\endref

\r\no \refGarside 
\by     F.~A.~Garside 
\paper  The braid group and other groups 
\jour   Quart. J. Math. \vol 20 \yr 1969 \pages 235--254 
\endref 

\r\no\refGebhardt
\by V.~Gebhardt
\paper A new approach to the conjugacy problem in Garside groups
\jour J. of Algebra \vol 292 \yr 2005 \pages 282--302
\endref

\r\no\refGebGM
\by V.~Gebhardt, J.~Gonz\'alez-Meneses
\paper The cyclic sliding operation in Garside groups
\jour  Math. Z. \vol 265 \yr 2010 \pages 85--114
\endref

\r\no\refMurasugi
\by K.~Murasugi
\book On closed 3-braids \bookinfo Memoirs of the AMS, 151
\publ Amer. Math. Soc. \publaddr Providence, RA \yr 1974
\endref

\r\no\refOrevkovTop
\by S.~Yu.~Orevkov
\paper Link theory and oval arrangements of real algebraic curves
\jour Topology \vol 38 \yr 1999 \pages 779--810
\endref

\r\no\refOrevkovUR
\by S.~Yu.~Orevkov
\paper Quasipositivity test via unitary representations of braid
groups and its applications to real algebraic curves
\jour J. Knot Theory and Ramifications \vol 10 \yr 2001 \pages 1005--1023
\endref

\r\no\refOrevkovGAFA
\by S.~Yu.~Orevkov
\paper Classification of flexible $M$-curves of degree $8$ up to isotopy
\jour GAFA - Geom. and Funct. Anal. \vol 12 \yr 2002 \pages 723--755
\endref

\r\no\refOrevkovQPthree
\by S.~Yu.~Orevkov
\paper Quasipositivity problem for 3-braids
\jour Turkish J. Math. \vol 28 \yr 2004 \pages 89--93
\endref

\r\no\refOrevkovDN
\by S.~Yu.~Orevkov
\paper Plane real algebraic curves of odd degree with a deep nest 
\jour J. Knot Theory and Ramifications \vol 14 \yr 2005 \pages 497--522
\endref

\r\no\refOrevkovAa
\by S.~Yu.~Orevkov
\paper Arrangements of an $M$-quintic with respect to a conic
       which maximally intersects its odd branch
\jour  Algebra i Analiz \vol 19:4 \yr 2007 \pages 174--242 \lang Russian
\transl English transl. St. Petersburg Math. J. \vol 19 \yr 2008 \pages 625--674
\endref

\r\no\refOrW
\by S.~Yu.~Orevkov
\paper Some examples of real algebraic and real pseudoholomorphic curves
\inbook in: Perspectives in Analysis, Geometry and Topology
\bookinfo Progr. in Math. 296 \publ Birkh\"auser/Springer \publaddr N. Y.
\yr 2012 \pages 355-387 
\endref

\r\no\refPrasolov
\by M.~V.~Prasolov
\paper Small braids with large ultra summit set
\jour Mat. Zametki \vol 89 \yr 2011 \pages 577--588 \lang Russian
\transl English transl. \jour Math. Notes \vol 89:4 \yr 2011 \pages 545--554
\endref

\r\no\refRudolph
\by L.~Rudolph
\paper Algebraic functions and closed braids
\jour Topology \vol 22 \yr 1983 \pages 191--202
\endref

\endRefs

\enddocument